\documentclass[12pt,reqno]{amsart}
\usepackage[utf8x]{inputenc}
\usepackage{amssymb,amsmath,latexsym, tikz}
\usepackage{epsfig}

\usepackage{epic,eepic,wrapfig,color}

\newtheorem{thm}{Theorem}[section]
\newtheorem{cor}[thm]{Corollary}
\newtheorem{lemma}[thm]{Lemma}
\newtheorem{prop}[thm]{Proposition}
\newtheorem{defn}[thm]{Definition}

\newtheorem{remark}[thm]{Remark}
\newtheorem{example}[thm]{Example}
\numberwithin{equation}{section}

\def\pf{{\medskip\noindent {\bf Proof. }}}
\def\qed{{\hfill $\Box$ \bigskip}}
\def\R{{\mathbb R}}
\def\P{{\mathbb P}}
\def\E{{\mathbb E}}

\def\1{{\bf 1}}

\DeclareMathOperator{\dist}{dist}

\newcommand{\cal}[1]{\mathcal{#1}}

  \def\sL {{\cal L}}

  \def\bR {{\mathbb R}}

\def\R {{\mathbb R}}

\def\nn{\nonumber}

\def\E{{\mathbb E}}
\def\P{{\mathbb P}}

\def\bea{\begin{align*}}
\def\eea{\end{align*}}
\def\bee{\begin{equation}}
\def\eee{\end{equation}}

\makeatletter
\@addtoreset{equation}{section}

\makeatother

\begin{document}
\bibliographystyle{plain}

\title[
On estimates of Dirichlet heat kernel  for SBM] {\bf
On estimates of transition density for subordinate Brownian motions with Gaussian components in $C^{1,1}$-open sets}

\author{Joohak Bae}
\address{Department of Mathematical Sciences,
Seoul National University,
Building 27, 1 Gwanak-ro, Gwanak-gu
Seoul 08826, Republic of Korea}
\curraddr{}
\thanks{
This work was supported by the National Research Foundation of Korea(NRF) grant funded by the Korea government(MSIP) 
(No. NRF-2015R1A4A1041675).
}

\author{Panki Kim}
\address{Department of Mathematical Sciences,
Seoul National University,
Building 27, 1 Gwanak-ro, Gwanak-gu
Seoul 08826, Republic of Korea}
\curraddr{}
\thanks{
This work was supported by the National Research Foundation of Korea(NRF) grant funded by the Korea government(MSIP) (No. 2016R1E1A1A01941893).
}
\email{pkim@snu.ac.kr}

\email{juhak88@snu.ac.kr}


\maketitle

\begin{abstract}
We consider a subordinate Brownian motion $X$ with Gaussian components when the scaling order of purely discontinuous part is between $0$ and $2$ including $2$. In this paper we establish sharp two-sided bounds for transition density of $X$ in $\R^d$ and $C^{1,1}$-open sets. As a corollary, we obtain a sharp Green function estimates.
\end{abstract}



\bigskip

\section{Introduction}
One of the most important notions in probability theory and analysis is the heat kernel. The transition density $p(t,x,y)$ of a  Markov process $X$ is the heat kernel of infinitesimal generator $\sL$ of $X$(also called the fundamental solution of $\partial_tu=\sL u$), whose explicit form does not exists usually. Thus obtaining sharp estimates of heat kernel $p(t,x,y)$ is a fundamental problem in both fields. Recently, for a large class of purely discontinuous Markov processes, the sharp heat kernel estimates were obtained in \cite{CK2, CK3, CKK, CKK2}. A common property of all purely discontinuous Markov processes considered so far in the estimates of the heat kernel was that the scaling order was always strictly between 0 and 2. In \cite{M}, Ante Mimica succeeded in obtaining sharp heat kernel estimates for purely discontinuous subordinate Brownian motions when the scaling order is between 0 and 2 including 2. For heat kernel estimates of processes with diffusion parts, mixture of Brownian motion and stable process was considered in \cite{SV} and diffusion process with jumps was considered in \cite{CK}.

 For any open subset $D\subset \bR^d$, let  $X^D$ be a subprocess of $X$ killed upon leaving $D$ and $p_D(t, x, y)$ be a transition density of $X^D$.  
An infinitesimal generator $\sL|_D$ of $X^D$ is the infinitesimal generator $\sL$ with zero exterior condition. 
$p_D(t,x,y)$ is  also called the Dirichlet heat kernel for $\sL|_D$ since 
it is the fundamental solution to  exterior Dirichlet problem with respect to $\sL|_D$. There are many results for Dirichlet heat kernel estimates in open subsets of $\R^d$ (see \cite{CKS7, CKS2, CKS3, CKS4, CKS5, CKS6}). In \cite{CKS7}, second-named author, jointly with Zhen-Qing Chen and Renming Song, obtains sharp two-sided estimates for the Dirichlet heat kernels in $C^{1,1}$-open sets of a large class of subordinate Brownian motions with Gaussian components.
Very recently, In \cite{KM} second-named author, jointly with Ante Mimica, establish sharp two-sided estimates
 for the Dirichlet heat kernels in $C^{1,1}$-open sets of  subordinate Brownian motions without Gaussian components whose scaling order is not necessarily strictly below $2$.

 In this paper, we continue the journey  on investigating the sharp two-sided estimates of heat kernels both in the whole space and $C^{1,1}$-open sets.
Here, we consider subordinate Brownian motions with Gaussian components when the scaling order of purely discontinuous part is between 0 and 2 including 2. Such processes were not considered in \cite{SV, CK, CKS7, KM}.

Let us describe the results of the paper in more detail. We start with a
description of the setup of this paper.

Let $S=(S_t)_{t\geq0}$ be a subordinator (increasing 1-dimensional L\'evy process) whose Laplace transform of $S_t$ is of the form
$$
\E e^{-\lambda S_t}=e^{-t\psi(\lambda)}, \quad \lambda>0,
$$
where $\psi$ is called the Laplace exponent of $S$. Without loss of generality, we assume the drift of $\psi$ is equal to 1 so that $\psi$ has the expression 
\begin{equation}\label{lapex}
\psi(\lambda)=\lambda + \phi(\lambda)\quad \text{with}\quad \phi(\lambda):=\int_{(0,\infty)}(1-e^{-\lambda t})\mu(dt).
\end{equation}
Here, $\mu$ is a L\'evy measure of $S$ satisfying $\int_{(0,\infty)}(1\wedge t)\mu(dt)<\infty$.\\

Let $X=(X_t)_{t\geq0}=(W_{S_t})_{t\geq 0}$ be a subordinate Brownian motion with subordinator $S=(S_t)_{t\geq0}$, where $W=(W_t)_{t\geq 0}$ is a Brownian motion independent of $S$. Then $X$ is rotationally invariant L\'evy process whose characteristic function is $\psi(|\xi|^2)=|\xi|^2+\phi(|\xi|^2)$.
One can view $X$ as an independent sum of a Brownian motion and purely discontinuous subordinate Brownian motion i.e., $X_t=B_t+Y_t$ where $B$ is a Brownian motion and $Y$ is a subordinate Brownian motion, independent of $B$, with subordinator $T$ whose Laplace exponent of $T$ is $\phi$. If the scaling order of $\phi$ is 2,  one can say that the process $X$ is very close to Brownian motion. (See Corollary \ref{ge}.)\\
The L\'evy density (jumping kernel) $J$ of $X$ is given by
$$
J(x)=j(|x|)=\int_0^{\infty}(4\pi t)^{-d/2}e^{-|x|^2/4t}\mu(t)dt.
$$
The function $J(x)$ determines a L\'evy system for $X$: for any non-negative measurable function $f$ on $\R_+ \times \R^d \times \R^d$ with $f(s,y,y)=0$ for all $y\in \R^d$, any stopping time $T$ (with respect to the filtration of $X$) and any $x\in \R^d$,
\begin{equation}\label{levy}
\E_x \left[\sum_{s\leq T}f(s,X_{s-},X_s)\right]=\E_x\left[\int_0^T\left(\int_{\R^d}f(s,X_s,y)J(X_s-y)dy\right)ds\right].
\end{equation}
\\
We first introduce the following scaling conditions for a function $f:(0,\infty) \to (0,\infty)$.
\begin{defn}{\rm
Suppose $f$ is a function from $(0,\infty)$ to $(0,\infty)$.\\
(1) We say that $f$ satisfies $L_a(\gamma, C_L)$ if there exist $a\geq 0$, $\gamma>0$, and $C_L\in (0,1]$ such that
$$ \frac{f(\lambda x)}{f(\lambda)} \geq C_L x^\gamma \quad \text{for all} \quad \lambda > a \quad \text{and} \quad x\geq 1,$$
(2) We say that $f$ satisfies $U_a(\delta, C_U)$ if there exist $a\geq 0$, $\delta>0$, and $C_U\in [1,\infty)$ such that
$$ \frac{f(\lambda x)}{f(\lambda)} \leq C_U x^\delta \quad \text{for all} \quad \lambda > a \quad \text{and} \quad x\geq 1.$$}
\end{defn}
\begin{remark}\label{msc}
According to \cite[Remark 2.2]{KM}, if we assume in addition $f$ is increasing, then the following holds.\\
(1) If $f$ satisfies $L_b(\gamma, C_L)$ with $b>0$ then $f$ satisfies $L_a(\gamma, (\frac{a}{b})^{\gamma}C_L)$ for all $a\in (0,b]$:
$$
\frac{f(\lambda x)}{f(\lambda)} \geq \Big(\frac{a}{b}\Big)^{\gamma} C_L x^{\gamma}, \quad x\geq1, \lambda \geq a.
$$
(2) If $f$ satisfies $U_b(\delta, C_U)$ with $b>0$ then $f$ satisfies $U_a(\delta, \frac{f(b)}{f(a)}C_U)$ for all $a\in (0,b]$:
$$
\frac{f(\lambda x)}{f(\lambda)}\leq \frac{f(b)}{f(a)}C_U x^{\delta}, \quad x\geq1, \lambda \geq a.
$$
\end{remark}
Throughout this paper we denote $p^{(2)}(t,x)$ the transition density of $B$ (and $W$).
i.e., $$
p^{(2)}(t,x)=(4\pi t)^{-d/2}\exp(-\frac{|x|^2}{4t}).
$$
We assume that $\mu(0, \infty)=\infty$ and denote $q(t,x)$ the transition density of $Y$ and $p(t,x)$ the transition density of $X$.
$q(t,x)$  and $p(t,x)$ are  of the forms
\begin{align}\label{eq:heat_kernel}
p(t,x)&=\int_{(0,\infty)}(4\pi s)^{-d/2}e^{-\frac{|x|^2}{4s}}\P(S_t\in ds), \quad
q(t,x)&=\int_{(0,\infty)}(4\pi s)^{-d/2}e^{-\frac{|x|^2}{4s}}\P(T_t\in ds)
\end{align}
for $x,y\in \R^d$ and $t>0$\,. These imply that for all $t>0$, $p(t,x) \le p(t,y)$ and $q(t,x) \le q(t,y)$  if $|x| \ge |y|$.\\
Throughout this paper the constants $r_0, R_0, \lambda_0, \Lambda_0$, and $C_i, i=1,2,...$ will be fixed. While, we use $c_1, c_2, ...$ to denote generic constants, whose exact values are not important and the labeling of the constants $c_1, c_2,...$ starts anew in the statement of each result and its proof. For $a,b \in \R$ we denote $a\wedge b:=\min\{a,b\}$ and $a\vee b:=\max\{a,b\}$. Notation $f(x)\asymp g(x), x\in I$ means that there exists constants $c_1, c_2>0$ such that $c_1 f(x)\leq g(x) \leq c_2 g(x)$ for $x\in I$.
\\
\\
The following is the first main result of this paper. Throughout this paper, 
$S=(S_t)_{t\geq0}$ is a subordinator whose Laplace exponent $\psi$ is $\lambda + \phi(\lambda)$ 
and we denote $H(\lambda):=\phi(\lambda)-\lambda\phi'(\lambda)$. 
\begin{thm}\label{T:2}
Let $X=(X_t)_{t\geq0}$ be a subordinate Brownian motion whose characteristic exponent is $\psi(|\xi|^2)=|\xi|^2+\phi(|\xi|^2)$.\\
(1) Suppose $H$ satisfies $L_a(\gamma, C_L)$ and $U_a(\delta, C_U)$ with $\delta <2$ for some $a>0$. For every $T,M>0$, there exist positive constants $c_1, c_2, c_3, c_4$, and $c_5$ such that for all $0<t\leq T$ and $|x|\leq M/2$,
\begin{align*}
c_1^{-1}\big(t^{-d/2}\wedge \big(p^{(2)}(t,c_2x)+q(t,c_3x)\big)\big)\leq
p(t,x)
\leq c_1\big(t^{-d/2}\wedge \big(p^{(2)}(t,c_4x)+q(t,c_5x)\big)\big).
\end{align*}
(2) Suppose $H$ satisfies $L_0(\gamma, C_L)$ and $U_0(\delta, C_U)$ with $\delta <2$. Then there exist positive  constants $c_6, c_7, c_8, c_9$, and $c_{10}$ such that for all $t>0$ and $x\in \R^d$,
\begin{align*}
c_6^{-1}\Big(\big(t^{-d/2}\wedge & \phi^{-1}(t^{-1})^{d/2}\big)\wedge \big(p^{(2)}(t,c_7x)+q(t,c_8x)\big)\Big)\leq
p(t,x)\\
&\leq c_6\Big(\big(t^{-d/2}\wedge \phi^{-1}(t^{-1})^{d/2}\big)\wedge \big(p^{(2)}(t,c_9x)+q(t,c_{10}x)\big)\Big).
\end{align*}
\end{thm}
As an application, we obtain sharp two-sided estimate for Green function of transient subordinate Brownian motion $X=(X_t)_{t\geq 0}$ $(d\geq 3)$. If $X$ is transient, then the following Green function is well-defined and finite.
$$
G(x,y)=G(x-y)=\int_0^{\infty}p(t,x-y)dt, \quad x,y\in \R^d, x\not=y.
$$
\begin{cor}\label{ge}
Let $d\geq3$. Suppose $H$ satisfies $L_0(\gamma, C_L)$ and $U_0(\delta, C_U)$ with $\delta <2$. Then
$$
G(x) \asymp |x|^{-d} \left(|x|^{2}\wedge \phi(|x|^{-2})^{-1}\right), \quad x\in \R^d.
$$
\end{cor}
Consider $\phi$ is given in Example \ref{Ex12} (2) below. Then
$$
|x|^{-d}\phi(|x|^{-2})^{-1} \asymp \begin{cases}
|x|^{-d+2}\log\frac{1}{|x|} \quad &|x|<\frac{1}{2} \\
|x|^{-d+2} &|x|\geq \frac{1}{2}.
\end{cases}
$$
Thus $G(x)\asymp G^{(2)}(x)$, $x\in \R^d$ where $G^{(2)}(x)=c|x|^{-d+2}$ is the Green function of the Brownian motion.
This shows that how close this process is to the Brownian motion and Green function estimates may not detect the difference between our $X$ and the Brownian motion.

Let $D\subset \R^d$ (when $d\geq 2$) be a $C^{1,1}$ open set with $C^{1,1}$ characteristics $(R_0,\Lambda_0)$, that is, there exists a localization radius $R_0 >0$ and a constant $\Lambda_0 >0$ such that for every $z \in \partial D$ there exist a $C^{1,1}$-function $\varphi =\varphi_z:\R^{d-1} \to \R$ satisfying $\varphi(0)=0$, $\nabla\varphi(0)=(0,...,0)$, $||\nabla \varphi||_{\infty}\leq \Lambda_0$, $|\nabla\varphi(x)-\nabla\varphi(w)|\leq \Lambda_0|x-w|$ and an orthonormal coordinate system $CS_z$ of $z=(z_1,\cdots,z_{d-1},z_d):=(\tilde{z},z_d)$ with origin at $z$ such that $D\cap B(z,R_0)=\{y=(\tilde{y},y_d) \in B(0,R_0)$ in  $CS_z:y_d>\varphi(\tilde{y})\}$. The pair $(R_0, \Lambda_0)$ will be called the $C^{1,1}$ characteristics of the open set $D$.
By a $C^{1,1}$ open set  in $\R$ with a characteristic $R_0>0$, we mean an open set that can be written as the union of disjoint intervals so that the {infimum} of the lengths of all these intervals is {at least $R_0$} and the {infimum} of the distances between these intervals is  {at least $R_0$}.

Throughout this paper we  denote  $p_D(t,x,y)$ the transition density of $X^D$.
The following are the second main results in this paper.

\begin{thm}\label{T:3}
Let $X=(X_t)_{t\geq0}$ be a subordinate Brownian motion whose characteristic exponent is $\psi(|\xi|^2)=|\xi|^2+\phi(|\xi|^2)$.
Suppose $H$ satisfies $L_a(\gamma, C_L)$ and $U_a(\delta, C_U)$ with $\delta <2$ for some $a>0$ and $D$ is a bounded $C^{1,1}$ open set in $\R^d$ with characteristics $(R_0,\Lambda_0)$. Then for every $T>0$ there exist positive constants $c_1, c_2, c_3, a_U, a_L$ such that \\
(1) For any $(t,x,y)\in (0,T]\times D \times D$, we have
\begin{align}\label{dhkeu}
&p_D(t,x,y) \leq c_1\bigg(1 \wedge \frac{\delta_D(x)}{\sqrt{t}}\bigg)\bigg(1 \wedge \frac{\delta_D(y)}{\sqrt{t}}\bigg)\nonumber  \\
&\quad \times \Big(t^{-d/2}\wedge (p^{(2)}(t,c_2(x-y)) + \frac{tH(|x-y|^{-2})}{|x-y|^{d}}+\phi^{-1}(t^{-1})^{d/2}e^{-a_U|x-y|^2\phi^{-1}(t^{-1})})\Big).
\end{align}
(2) For any $(t,x,y)\in(0,T] \times D \times D$, we have
\begin{align}\label{dhkel}
&p_D(t,x,y) \geq c_1^{-1}\bigg(1 \wedge \frac{\delta_D(x)}{\sqrt{t}}\bigg)\bigg(1 \wedge \frac{\delta_D(y)}{\sqrt{t}}\bigg)\nonumber\\
&\quad \times \Big(t^{-d/2}\wedge (p^{(2)}(t,c_3(x-y))+\frac{tH(|x-y|^{-2})}{|x-y|^{d}}+\phi^{-1}(t^{-1})^{d/2}e^{-a_L|x-y|^2\phi^{-1}(t^{-1})})\Big).
\end{align}
(3) For any $(t,x,y)\in [T,\infty)\times D\times D$, we have
$$
p_D(t,x,y) \asymp e^{-\lambda_1t}\delta_D(x)\delta_D(y),
$$
where $-\lambda_1<0$ is the largest eigenvalue of the generator of $X^D$.
\end{thm}

We say that {\it the path distance in a domain (connected open set) $U$ is comparable to the Euclidean distance with characteristic $\lambda_0$} if for every $x$ and $y$ in  $U$ there is a rectifiable curve $l$ in $U$ which connects $x$ to $y$ such that the length of $l$ is less than or equal to  $\lambda_0|x-y|$.
Clearly, such a property holds for all bounded $C^{1,1}$ domains, $C^{1,1}$ domains with compact complements, and domain consisting of all the points above the graph of a 
bounded globally $C^{1,1}$ function. 

\begin{thm}\label{T:4}
Let $X=(X_t)_{t\geq0}$ be a subordinate Brownian motion whose characteristic exponent is $\psi(|\xi|^2)=|\xi|^2+\phi(|\xi|^2)$. Suppose $H$ satisfies $L_0(\gamma, C_L)$ and $U_0(\delta, C_U)$ with $\delta <2$ and $D$ is an unbounded $C^{1,1}$ open set in $\R^d$ with characteristics $(R_0,\Lambda_0)$. Then for every $T>0$ there exists $c_1, c_2, c_3, a_U, a_L$ such that \\
(1) For any $(t,x,y)\in (0,T]\times D \times D$, we have
\begin{align*}
&p_D(t,x,y) \leq c_1\bigg(1 \wedge \frac{\delta_D(x)}{\sqrt{t}}\bigg)\bigg(1 \wedge \frac{\delta_D(y)}{\sqrt{t}}\bigg) \\
&\quad\quad\quad\quad\quad \times  \Big(t^{-d/2}\wedge (p^{(2)}(t,c_2(x-y)) + \frac{tH(|x-y|^{-2})}{|x-y|^d}+\phi^{-1}(t^{-1})^{d/2}e^{-a_U|x-y|^2\phi^{-1}(t^{-1})})\Big).
\end{align*}
(2) If the path distance in each connected component of $D$ is comparable to the Euclidean distance with characteristic $\lambda_0$, then for any $(t,x,y)\in(0,T] \times D \times D$, we have
\begin{align*}
&p_D(t,x,y) \geq c_1^{-1}\bigg(1\wedge \frac{\delta_D(x)}{\sqrt{t}}\bigg)\bigg(1\wedge \frac{\delta_D(y)}{\sqrt{t}}\bigg)\\
&\quad\quad\quad\quad\quad \times \Big(t^{-d/2}\wedge (p^{(2)}(t,c_3(x-y))+\frac{tH(|x-y|^{-2})}{|x-y|^d}+\phi^{-1}(t^{-1})^{d/2}e^{-a_L|x-y|^2\phi^{-1}(t^{-1})})\Big).
\end{align*}

\end{thm}
Define $G_D(x,y)=\int_0^{\infty}p_D(t,x,y)dt$, Green function of $X^D$. The following is Green function estimate of $X^D$.

\begin{cor}\label{dge}
Suppose $H$ satisfies $L_a(\gamma, C_L)$ and $U_a(\delta, C_U)$ with $\delta <2$ for some $a>0$ and $D$ is a bounded $C^{1,1}$ open set in $\R^d$ with characteristics $(R_0,\Lambda_0)$. Then
$$
G_D(x,y) \asymp g_D(x,y), \quad x,y\in D,
$$
where 
\begin{align}\label{g_D}
g_D(x,y):= \begin{cases}
\frac{1}{|x-y|^{d-2}}\bigg(1\wedge \frac{\delta_D(x)\delta_D(y)}{|x-y|^2}\bigg) \quad &\text{when} \quad d\geq 3,\\
\log\bigg(1+\frac{\delta_D(x)\delta_D(y)}{|x-y|^2}\bigg) &\text{when} \quad d=2,\\
\big(\delta_D(x)\delta_D(y)\big)^{1/2}\wedge \frac{\delta_D(x)\delta_D(y)}{|x-y|} &\text{when}\quad d=1.
\end{cases}
\end{align}
\end{cor}
Denote by $G_D^{(2)}(x,y)$ the Green function of Brownian motion in $D$. It is known (see \cite{CZ}) that $G_D^{(2)} \asymp g_D(x,y)$ when $x$ and $y$ are in the same component of $D$, and $G_D^{(2)}(x,y)=0$ otherwise. Thus when $D$ is a bounded $C^{1,1}$ connected open subset of $\R^d$, $G_D(x,y) \asymp G_D^{(2)}(x,y)$, while our heat kernel estimates (Theorem \ref{T:3}) are different from heat kernel estimates of Brownian motion in $D$.\\
\\
These are examples where the scaling order of $\phi$ is not strictly between $0$ and $2$.
\begin{example}\label{Ex12}
(1) Let $\phi(\lambda)=\frac{\lambda}{\log(1+\lambda^{\beta/2})}$, where $\beta \in (0,2)$. Then
\begin{align*}
\phi^{-1}(\lambda) \asymp 
\left\{ \begin{array}{ll} \lambda^{\frac{2}{2-\beta}} \quad &0<\lambda<2\\
\lambda \log\lambda \quad &\lambda \geq 2\\
\end{array} \right. \quad \quad
H(\lambda)\asymp
\left\{ \begin{array}{ll} \lambda^{1-\beta/2} \quad &0<\lambda<2\\
\frac{\lambda}{(\log\lambda)^2} \quad &\lambda \geq 2\\
\end{array} \right.
\end{align*}
Hence, $H$ satisfies $L_0(\gamma, C_L)$ and $U_0(\delta, C_U)$ with some $\gamma, C_L, C_U$ and $\delta<2$.\\
\\
(2) Let $\phi(\lambda)=\frac{\lambda}{\log(1+\lambda)}-1$. Then
\begin{align*}
\phi^{-1}(\lambda) \asymp 
\left\{ \begin{array}{ll} \lambda \quad &0<\lambda<2\\
\lambda \log\lambda \quad &\lambda \geq 2\\
\end{array} \right. \quad \quad
H(\lambda)\asymp
\left\{ \begin{array}{ll} \lambda^2 \quad &0<\lambda<2\\
\frac{\lambda}{(\log\lambda)^2} \quad &\lambda \geq 2\\
\end{array} \right.
\end{align*}
Hence, $H$ satisfies $L_0(\gamma, C_L)$ and $U_2(\delta, C_U)$ with some $\gamma, C_L, C_U$ and $\delta<2$.\\
\\
 Suppose $D$ is a bounded $C^{1,1}$ open set with $diam(D)<1/2$ and $\psi(\lambda)=\lambda + \phi(\lambda)$, where $\phi$ is the one in  above two cases. Then for $t<1/2$, there exist positive constants $c_1, c_2, c_3, a_U$, and $a_L$ such that
\begin{align*}
&p_D(t,x,y) \leq c_1\bigg(1 \wedge \frac{\delta_D(x)}{\sqrt{t}}\bigg)\bigg(1 \wedge \frac{\delta_D(y)}{\sqrt{t}}\bigg) \\
&\quad\quad \times \left(t^{-d/2}\wedge \Big(p^{(2)}(t,c_2(x-y)) + \frac{t}{|x-y|^{d+2}(\log\frac{1}{|x-y|})^2}+t^{-d/2}\Big(\log\frac{1}{t}\Big)^{d/2}e^{-a_U\frac{|x-y|^2}{t}\log\frac{1}{t}}\Big)\right),
\end{align*}
and
\begin{align*}
&p_D(t,x,y) \geq c_1^{-1}\bigg(1 \wedge \frac{\delta_D(x)}{\sqrt{t}}\bigg)\bigg(1 \wedge \frac{\delta_D(y)}{\sqrt{t}}\bigg) \\
&\quad\quad \times \left(t^{-d/2}\wedge \Big(p^{(2)}(t,c_3(x-y)) + \frac{t}{|x-y|^{d+2}(\log\frac{1}{|x-y|})^2}+t^{-d/2}\Big(\log\frac{1}{t}\Big)^{d/2}e^{-a_L\frac{|x-y|^2}{t}\log\frac{1}{t}}\Big)\right).
\end{align*}

\end{example}

\section{Heat kernel estimates in $\R^d$}\label{S:int}
Throughout this paper, 
$X=(X_t)_{t\geq0}$ is a subordinate Brownian motion whose characteristic exponent is $\psi(|\xi|^2)=|\xi|^2+\phi(|\xi|^2)$.
In this section we obtain estimates of transition density of the subordinate Brownian motion $X$.
The following are heat kernel estimates for $q(t,x)$, which is transition density of $Y$. 
Recall that $Y$ is a subordinate Brownian motion with subordinator $T$ whose Laplace exponent of $T$ is $\phi$ and $H(\lambda)=\phi(\lambda)-\lambda\phi'(\lambda)$.

\begin{thm}[\cite{M, KM}]\label{t:MKM}
{\rm (i)} If $\phi$ satisfies $L_a(\gamma, C_L)$ for some $a>0$, then for every $T>0$ there exist $C_1=C_1(T)>1$ and $a_U>0$ such that for all $t\leq T$ and $x\in \R^d$, 
\begin{equation}\label{ouhke}
q(t,x)\leq C_1 \left(\phi^{-1}(t^{-1})^{d/2}\wedge \big(t|x|^{-d}H(|x|^{-2})+ \phi^{-1}(t^{-1})^{d/2}e^{-a_U|x|^2\phi^{-1}(t^{-1})}\big) \right),
\end{equation}
and
\begin{equation}\label{nhke}
 q(t,x) \ge C_1^{-1}\phi^{-1}(t^{-1})^{d/2}, \quad \text{if} \quad t\phi(|x|^{-2})\geq 1.
\end{equation}
Consequently, the L\'evy density (jumping kernel) $J$ satisfies 
$$J(x)=\lim_{t\to 0}q(t,x)/t\leq C_1|x|^{-d}H(|x|^{-2}).$$  Furthermore, if $a=0$, then (\ref{ouhke}) and (\ref{nhke}) hold for every $t>0$ and $x\in \R^d$. \\
{\rm (ii)} If $H$ satisfies $L_a(\gamma, C_L)$ and $U_a(\delta, C_U)$ with $\delta <2$ for some $a>0$, then for every $T, M>0$ there exist $C_1=C_1(a,\gamma, C_L, \delta, C_U, T, M)>0$ and $a_L>0$ such that for all $t\leq T$ and $|x|<M$,\begin{equation}\label{olhke}
q(t,x) \geq C_1^{-1} \left(\phi^{-1}(t^{-1})^{d/2}\wedge \big(t|x|^{-d}H(|x|^{-2})+ \phi^{-1}(t^{-1})^{d/2}e^{-a_L|x|^2\phi^{-1}(t^{-1})}\big) \right).
\end{equation}
Consequently, the L\'evy density $J$ satisfies \begin{equation}\label{jl}J(x)=\lim_{t\to 0}q(t,x)/t\asymp |x|^{-d}H(|x|^{-2}), \quad |x|<M.
\end{equation}
Furthermore,  if $a=0$, then (\ref{olhke}) and (\ref{jl}) hold for all $t>0$ and $x\in \R^d$. 
\end{thm}
We will use following formula of transition density $p(t,x)$ of $X$, which is given by
$$
p(t,x)=\int_{\R^d}p^{(2)}(t,x-y)q(t,y)dy.
$$
\\

\subsection{Upper bounds.} In this subsection we will prove the upper bounds for the transition density. First, we observe the following simple upper bound of $p(t,x)$. (See \cite{SV}.)
\begin{lemma}\label{L:2.4}
For every $(t, x) \in (0, \infty) \times \R^d$, $$
p(t,x) \leq \exp(|x|^2/(4t)) p^{(2)}(t,x).
$$
In particular if $t\geq |x|^2$, we have $p(t,x)\leq e^{1/4}p^{(2)}(t,x)$.
\end{lemma}

\pf
Since $p^{(2)}(t,x-y) \le  p^{(2)}(t,x) \exp (|x|^2/(4t)),$ we have 
\begin{align*}
p(t,x)
& = \int_{\R^d}p^{(2)}(t,x-y)q(t,y)dy \\
& \leq p^{(2)}(t,x)\exp(|x|^2/(4t))\int_{\R^d}q(t,y)dy =\exp(|x|^2/(4t)) p^{(2)}(t,x). 
\end{align*}
\qed\\
The next two lemmas will be used several times in this paper. 
\begin{lemma}[{\cite[Lemma 2.1(a)]{M}}]\label{lem:bf}
For any $\lambda>0$ and $x\geq 1$,
\[
    \phi(\lambda x)\leq x\phi(\lambda)\qquad \text{ and }\qquad H(\lambda x)\leq x^2H(\lambda)\,.
\]
\end{lemma}

\begin{lemma}[{\cite[Lemma 2.1(b)]{M}}]\label{lem:bf2}
For $a\geq0$ if $H$ satisfies $ L_a(\gamma, C_L)$ (resp. $ U_a(\delta, C_U)$), then $\phi$ satisfies $ L_a(\gamma, C_L)$(resp. $ U_a(\delta \wedge 1, C_U)$).  
\end{lemma}

Since the proofs are basically same, we provide the proof for the case $a>0$ in the next two results. 
\begin{lemma}\label{L:2.3}
Suppose $\phi$ satisfies $L_a(\gamma, C_L)$ for some $a>0$ ($L_0(\gamma, C_L)$, respectively). For $T>0$, there exists a positive constant $c$ such that for $0<t\leq T$ ($t>0$, respectively) and $t\phi(|x|^{-2})\geq 1$,
$$
p(t,x) \leq cq(t,x).
$$
\end{lemma}

\pf
Since $\phi$ satisfies $L_a(\gamma, C_L)$, for all $y \in \R^d$ and $0<t\leq T$, $q(t,y) \leq c_{1}\phi^{-1}(t^{-1})^{d/2}$ by Theorem \ref{t:MKM}(i). By (\ref{nhke}), for all $0<t\leq T$ and $t\phi(|x|^{-2})\geq 1$,  we have that $q(t,x)\geq C_{1}^{-1}\phi^{-1}(t^{-1})^{d/2}$. Hence,
\begin{align*}
p(t,x) 
& = \int_{\R^d}p^{(2)}(t,x-y)q(t,y)dy \\
& \leq c_1\phi^{-1}(t^{-1})^{d/2} \int_{\R^d}p^{(2)}(t,x-y)dy \leq c_{1} C_1q(t,x). 
\end{align*}
\qed
\\
Let 
$\tilde{p}^{(2)}(t,x):=(4\pi t)^{-d/2}\exp\{-|x|^2/(16t)\}.$
\begin{lemma} \label{L:2.5}
 Suppose $H$ satisfies $L_a(\gamma, C_L)$ and $U_a(\delta, C_U)$ for some $a>0$ 
($L_0(\gamma, C_L)$ and $U_0(\delta, C_U)$, respectively)
with $\delta <2$. For $T, M>0$ there exist positive constants $c$ and $c_0<1$ such that for all $0<t\leq T$ and $x\in \R^d$ satisfying $|x|<M$ ($t>0$ and $x\in \R^d$, respectively) and $t\phi(|x|^{-2})\leq 1$,
$$
p(t,x)\leq c\max (\tilde{p}^{(2)}(t,x), q(t,c_0x)).
$$
\end{lemma}
\pf 
We divide the integral
\begin{align*}
p(t,x)&=\int_{\R^d}p^{(2)}(t,x-y)q(t,y)dy\\
& = \int_{|y-x|>|x|/2}p^{(2)}(t,x-y)q(t,y)dy+\int_{|y-x|\leq |x|/2}p^{(2)}(t,x-y)q(t,y)dy \\
& =:I_1 + I_2.
\end{align*}
(i) For $|y-x|>|x|/2$, $\exp\{-|x-y|^2/(4t)\} \leq \exp\{-|x|^2/(16t)\}$. Therefore, 
\begin{align*}
I_1 &\leq (4\pi t)^{-d/2}\exp(-|x|^2/(16t))\int_{|y-x|>|x|/2}q(t,y)dy\leq \tilde{p}^{(2)}(t,x).
\end{align*}
(ii) For $|y-x|\leq |x|/2$, we have $\frac{|x|}{2} \leq |y|\leq \frac{3|x|}{2}$.
 Therefore using (\ref{ouhke}), Lemma \ref{lem:bf}, and (\ref{olhke}), we have
\begin{align*}
I_2 &\leq c_1\int_{|y-x|\leq|x|/2}p^{(2)}(t,x-y)(t|y|^{-d}H(|y|^{-2})+\phi^{-1}(t^{-1})^{d/2}e^{-a_U|y|^2\phi^{-1}(t^{-1})})dy\\
&\leq c_1\int_{|y-x|\leq|x|/2}p^{(2)}(t,x-y)(t2^d|x|^{-d}H(2^2|x|^{-2})+\phi^{-1}(t^{-1})^{d/2}e^{-a_U|x/2|^2\phi^{-1}(t^{-1})})dy\\
&\leq  c_2 \left(\int_{|y-x|\leq|x|/2}p^{(2)}(t,x-y)dy\right) \left(t|c_0x|^{-d}H(|c_0x|^{-2})+\phi^{-1}(t^{-1})^{d/2}e^{-a_L|c_0x|^2\phi^{-1}(t^{-1})}\right)\\
& \leq c_3q(t,c_0x),
\end{align*}
where $c_0=a^{1/2}_U/2a^{1/2}_L<1$.
\qed

\begin{remark}\label{R:L2.3}
In the proof of Lemma \ref{L:2.5}, we just used that $\phi$ satisfies $L_a(\gamma, C_L)$ until the last inequality. Thus if  $\phi$ satisfies $L_a(\gamma, C_L)$ for some $a>0$ ($L_0(\gamma, C_L)$, respectively), then for $T>0$ there exists a positive constant $c$ such that for all $0<t\leq T$ ($t>0$, respectively) and $x\in \R^d$ satisfying $t\phi(|x|^{-2})\leq 1$,
$$
p(t,x)\leq c\max(\tilde{p}^{(2)}(t,x), t|x|^{-d}H(|x|^{-2})+\phi^{-1}(t^{-1})^{d/2}e^{-(a_U/4)|x|^2\phi^{-1}(t^{-1})}).
$$
Consequently, combining Lemma \ref{L:2.4} and Lemma \ref{L:2.3} we can obtain upper bounds for $p(t,x)$:
If $\phi$ satisfies $L_a(\gamma, C_L)$ for some $a>0$, then for $0<t\leq T$ and $x \in \R^d$
$$
p(t,x)\leq c_1\Big(t^{-d/2}\wedge (t^{-d/2}e^{-|x|^2/(c_2t)} + t|x|^{-d}H(|x|^{-2})+\phi^{-1}(t^{-1})^{d/2}e^{-c_3|x|^2\phi^{-1}(t^{-1})})\Big),
$$
and if $\phi$ satisfies $L_0(\gamma, C_L)$ then for $t>0$ and $x \in \R^d$
$$
p(t,x)\leq c_4\Big((t^{-d/2}\wedge \phi^{-1}(t^{-1})^{d/2})\wedge (t^{-d/2}e^{-|x|^2/(c_5t)} + \frac{tH(|x|^{-2})}{|x|^{d}}+\phi^{-1}(t^{-1})^{d/2}e^{-a_6|x|^2\phi^{-1}(t^{-1})})\Big).
$$

\end{remark}
\subsection {Lower bounds.} In this subsection we will prove the lower bounds for the transition density. As the subsection 2.1, 
 we provide the proof for the case $a>0$ only.
 
Let $\hat{p}^{(2)}:=(4\pi t)^{-d/2}\exp(-|x|^2/t)$. 
\begin{lemma}\label{L:2.8}
Suppose $\phi$ satisfies $L_a(\gamma, C_L)$ for some $a>0$ ($L_0(\gamma, C_L)$, respectively). For $T>0$ there exists a positive constant $c$ such that for $0<t\leq T$ ($t>0$, respectively) and $x\in R^d$ satisfying $t\phi(|x|^{-2})\leq 1$,
$$
p(t,x)\geq c\hat{p}^{(2)}(t,x).
$$
\end{lemma}
\pf 
If $|y|\leq \phi^{-1}(t^{-1})^{-1/2}$, then $|y|\leq \phi^{-1}(t^{-1})^{-1/2}\leq |x|$. Therefore $|y-x|\leq 2|x|$ and hence $\exp(-|x-y|^2/(4t))\geq \exp(-|x|^2/t)$. By (\ref{nhke}), $q(t,y)\geq C_1^{-1}\phi^{-1}(t^{-1})^{d/2}$ for $0<t\leq T$. Thus, 
\begin{align*}
p(t,x) &\geq \int_{B(0,\phi^{-1}(t^{-1})^{-1/2})}p^{(2)}(t,x-y)q(t,y)dy \\
& \geq c_1(4\pi t)^{-d/2}\exp(-|x|^2/t)\phi^{-1}(t^{-1})^{d/2}(\phi^{-1}(t^{-1})^{-1/2})^d  =c_2\hat{p}^{(2)}(t,x). 
\end{align*}
\qed
\begin{lemma}\label{L:2.9}
Suppose $H$ satisfies $L_a(\gamma, C_L)$ and $U_a(\delta, C_U)$ with $\delta <2$ for some $a>0$ ($L_0(\gamma, C_L)$ and $U_0(\delta, C_U)$, respectively). For $T, M>0$ there exist positive constants $c$ and $\tilde{c_0}>1$ such that for all $0<t\leq T$, $x\in \R^d$ satisfying $|x|<\frac{2}{3}M$ ($t>0$ and $x \in \R^d$, respectively) and every $y\in B(x,|x|/2)$, it holds that 
$$
{q(t,y)}\geq  c q(t,\tilde{c_0}x).
$$
\end{lemma}
\pf 
Note that $y\in B(x,|x|/2)$ implies $|y|<\frac{3}{2}|x|<M$. We consider four cases separately.

\noindent
{\it Case (1):} When $t\phi(|x|^{-2})\geq 1$ and  $t\phi(|y|^{-2})\geq 1$.
 Using (\ref{nhke}),
$$q(t,x)\leq C_1\phi^{-1}(t^{-1})^{d/2}\leq C_1^2q(t,y).$$

\noindent
{\it Case (2):} When $t\phi(|x|^{-2})\geq 1$ and $t\phi(|y|^{-2})\leq 1$. Using (\ref{olhke}), $|y|<\frac{3}{2}|x|$, $|x|^2\phi^{-1}(t^{-1})\leq 1$, and (\ref{ouhke}), we have
\begin{align*}
q(t,y) &\geq C_1^{-1}\left(t|y|^{-d}H(|y|^{-2})+ \phi^{-1}(t^{-1})^{d/2}e^{-a_L |y|^2 \phi^{-1}(t^{-1})}\right) \geq C^{-1}\phi^{-1}(t^{-1})^{d/2}e^{-a_L|y|^2 \phi^{-1}(t^{-1})}\\
&\geq C_1^{-1}\phi^{-1}(t^{-1})^{d/2}e^{-\frac{9}{4}a_L|x|^2 \phi^{-1}(t^{-1})} \geq C_1^{-1}\phi^{-1}(t^{-1})^{d/2}e^{-\frac{9}{4}a_L}\geq c_1\phi^{-1}(t^{-1})^{d/2} \geq c_2q(t,x).
\end{align*}
{\it Case (3):} When $t\phi(|x|^{-2})\leq 1$ and $t\phi(|y|^{-2})\geq 1$. Using (\ref{nhke}), $tH(|x|^{-2})\leq t\phi(|x|^{-2})\leq 1$, $ |x|^{-d}\leq \phi^{-1}(t^{-1})^{d/2}$, and (\ref{ouhke}), we have
\begin{align*}
q(t,y) &\geq C_1^{-1}\phi^{-1}(t^{-1})^{d/2} \geq c_3(tH(|x|^{-2})|x|^{-d}+ \phi^{-1}(t^{-1})^{d/2}e^{-a_U|x|^2\phi^{-1}(t^{-1})}) \geq c_4q(t,x).
\end{align*}
\noindent
{\it Case (4):} When $t\phi(|x|^{-2})\leq 1$ and  $t\phi(|y|^{-2})\leq 1$. Using (\ref{olhke}), $|y|<3|x|/2$, and (\ref{ouhke}), we have
\begin{align*}
q(t,y) &\geq C_1^{-1}(t|y|^{-d}H(|y|^{-2})+ \phi^{-1}(t^{-1})^{d/2}e^{-a_L|y|^2 \phi^{-1}(t^{-1})})\\
&\geq C_1^{-1}(t(2/3)^d|x|^{-d}H(|\frac{3}{2}x|^{-2})+ \phi^{-1}(t^{-1})^{d/2}e^{-a_L |\frac{3}{2}x|^2 \phi^{-1}(t^{-1})})\\
&\geq c_5(t|\tilde{c_0}x|^{-d}H(|\tilde{c_0}x|^{-2})+ \phi^{-1}(t^{-1})^{d/2}e^{-a_U|\tilde{c_0}x|^2 \phi^{-1}(t^{-1})})\\
& \geq c_6q(t,\tilde{c_0}x),
\end{align*}
where $\tilde{c_0}=\frac{3}{2}(a_L/a_U)^{1/2}>1$.
Since $q(t,x)\geq q(t,\tilde{c_0}x)$ by (\ref{eq:heat_kernel}),  we conclude that
$q(t,y)\geq cq(t,\tilde{c_0}x) $ for some $c>0$.
\qed
\begin{lemma}\label{L:L2.10}
Suppose $H$ satisfies $L_a(\gamma, C_L)$ and $U_a(\delta, C_U)$ with $\delta <2$ for some $a>0$ ($L_0(\gamma, C_L)$ and $U_0(\delta, C_U)$, respectively). For $T, M>0$ there exist constants $c$ and $\tilde{c_0}>1$ such that for all $0<t\leq T$ and $x\in \R^d$ satisfying $|x|<\frac{2}{3}M$ ($t>0$ and $x \in \R^d$, respectively) and $t\leq |x|^2$, 
$$
p(t,x)\geq cq(t,\tilde{c_0}x).
$$
\end{lemma}
\pf 
By Lemma \ref{L:2.9}, $q(t,y)\geq c_1q(t,\tilde{c_0}x)$ for every $y\in B(x,|x|/2)$. Using this and change of variable,  we have
\begin{align*}
p(t,x) 
&= q(t,\tilde{c_0}x)\int_{\R^d}\frac{q(t,y)}{q(t,\tilde{c_0}x)}p^{(2)}(t,x-y)dy \\
&\geq c_1q(t,\tilde{c_0}x)\int_{B(x,|x|/2)}p^{(2)}(t,x-y)dy \\
&= c_1q(t,\tilde{c_0}x)\int_{B(0,t^{-1/2}\frac{|x|}{2})}p^{(2)}(1,u)du \\
&\geq c_1\left(\int_{B(0,1/2)}p^{(2)}(1,u)du \right)q(t,\tilde{c_0}x)  = c_2q(t,\tilde{c_0}x).
\end{align*}
\qed
\begin{lemma}\label{L:2.11}
Suppose $H$ satisfies $L_a(\gamma, C_L)$ and $U_a(\delta, C_U)$ with $\delta <2$ for some $a>0$ ($L_0(\gamma, C_L)$ and $U_0(\delta, C_U)$, respectively). For $T, M>0$ there exists a constant $c$ such that for all $1\leq t\leq T$ and $x\in \R^d$ satisfying $|x|<M/2$ ($t\geq 1$ and $x\in \R^d$, respectively) and $t\phi(|x|^{-2})\geq 1$,
$$
p(t,x)\geq cq(t,x).
$$
\end{lemma}
\pf 
Assume $t\phi(|x|^{-2})\geq 1$ and let $b=M\phi^{-1}(T^{-1})^{1/2}/2$. Note that 
we have $q(t,x)\leq C_1\phi^{-1}(t^{-1})^{d/2}$ by (\ref{ouhke}).

If $|y-x|\leq b\phi^{-1}(t^{-1})^{-1/2}$, then $|y|\leq |x-y|+|x| \leq (b+1)\phi^{-1}(t^{-1})^{-1/2}$ and $|y| \leq |x-y|+|x| \leq b\phi^{-1}(t^{-1})^{-1/2}+|x| \leq M$. 
Thus by (\ref{olhke}), we have 
$$q(t,y)\geq C_1^{-1}\phi^{-1}(t^{-1})^{d/2}e^{-a_L|y|^2\phi^{-1}(t^{-1})}\geq c_1\phi^{-1}(t^{-1})^{d/2}\quad \text {for } |y-x|\leq b\phi^{-1}(t^{-1})^{-1/2}.$$
Therefore, using the above inequality and change of variable
\begin{align*}
p(t,x) &\geq \int_{|x-y|\leq b\phi^{-1}(t^{-1})^{-1/2}}p^{(2)}(t,x-y)q(t,y)dy\\
&\geq c_2\phi^{-1}(t^{-1})^{d/2}\int_{|x-y|\leq b\phi^{-1}(t^{-1})^{-1/2}}p^{(2)}(t,x-y)dy\\
&=c_2\phi^{-1}(t^{-1})^{d/2}\int_{|u|\leq bt^{-1/2}\phi^{-1}(t^{-1})^{-1/2}}p^{(2)}(1,u)du\\
&\geq c_2\phi^{-1}(t^{-1})^{d/2}\int_{|u|\leq b\phi^{-1}(1)^{-1/2}}p^{(2)}(1,u)du\\
&\geq c_3\phi^{-1}(t^{-1})^{d/2}\\
&\geq c_4q(t,x).
\end{align*}
In the third inequality, we use   $\frac{\phi^{-1}(1)}{\phi^{-1}(t^{-1})}\geq t$ which follows from Lemma \ref{lem:bf}. 
  (See also \cite[Lemma 3.1(i)]{M}).
\qed	
\begin{lemma}\label{L:2.12}
Suppose $\phi$ satisfies $L_a(\gamma, C_L)$ for some $a\geq 0$. For $T\geq 1$ there exists a positive constant $c$ such that for all $t\leq T$ and $t\phi(|x|^{-2})\geq 1$,
$$
p(t,x)\geq cp^{(2)}(t,x).
$$
\end{lemma}
\pf 
We may assume that $a<\phi^{-1}(t^{-1})$ by Remark \ref{msc}. By the condition $L_a(\gamma, C_L)$ on $\phi$ (See also \cite[Lemma 3.1(ii)]{M}) and Lemma \ref{lem:bf}, we have for $t\leq T$,
\begin{align}\label{phisc}
C_L^{-1/\gamma}T^{1/\gamma}\phi^{-1}(t^{-1})\geq \phi^{-1}(Tt^{-1})\geq Tt^{-1}\phi^{-1}(1).
\end{align}
If $|y|\leq \phi^{-1}(t^{-1})^{-1/2}$, then $q(t,y)\geq C_1^{-1}\phi^{-1}(t^{-1})^{d/2}$ by (\ref{nhke}).
Also $|x-y|\leq |x|+|y|\leq 2\phi^{-1}(t^{-1})^{-1/2}$ and \eqref{phisc} imply 
$$
\exp(-\frac{|x-y|^2}{4t})\geq \exp(-\frac{4\phi^{-1}(t^{-1})^{-1}}{4t})\geq \exp{(-C_L^{-1/\gamma}\phi^{-1}(1)^{-1}T^{1/\gamma-1}}).
$$
Therefore for $t\leq T$ and $t\phi(|x|^{-2})\geq 1$,
\begin{align*}
p(t,x)&\geq \int_{|y|\leq \phi^{-1}(t^{-1})^{-1/2}}p^{(2)}(t,x-y)q(t,y)dy\\
&\geq c_1\phi^{-1}(t^{-1})^{d/2}\int_{|y|\leq \phi^{-1}(t^{-1})^{-1/2}}p^{(2)}(t,x-y)dy\\
&\geq c_1\phi^{-1}(t^{-1})^{d/2}(4\pi t)^{-d/2}e^{-C_L^{-1/\gamma}\phi^{-1}(1)^{-1}T^{1/\gamma-1}}\int_{|y|\leq \phi^{-1}(t^{-1})^{-1/2}}dy\\
&\geq c_2(4\pi t)^{-d/2}\geq c_2p^{(2)}(t,x).
\end{align*}
\qed\\
{\bf Proof of Theorem \ref{T:2}.}  Combining Lemma \ref{lem:bf2}, Lemma \ref{L:2.4}, Lemma \ref{L:2.3}, and Lemma \ref{L:2.5}, we get the upper bound of $p(t,x)$. Using Lemma \ref{lem:bf2}, Lemma \ref{L:2.8}, Lemma \ref{L:L2.10}, Lemma \ref{L:2.11}, and Lemma \ref{L:2.12}, we get the lower bound of $p(t,x)$. See Figure1.
 \qed
 
\begin{tikzpicture}
  \draw[->,line width=1pt] (0,0) -- (0,8) node[left] {$t$};
  \draw[->,line width=1pt] (0,0) -- (10,0) node[below] {$|x|$};
\draw[scale=1,samples=200,domain=0:7.3,smooth,variable=\x, dotted, line width=1pt] plot ({\x},{0.12*(\x)^2.1});
\draw[scale=1,samples=200,domain=0:10,smooth,variable=\x, line width=1pt] plot ({\x},{1.8*\x^0.5});
\draw[scale=1,samples=200,domain=0:5.3,smooth,variable=\x,dotted, line width=1pt] plot ({\x},{4.2});
\draw[scale=1,samples=200,domain=0:4.3,smooth,variable=\x,dotted, line width=1pt] plot ({5.4},{\x});
\node at (3,6) { $q(t,x)$};
\node at (8,6) { $q(t,x)$};
\node
 at (8,2) { $p^{(2)}(t,c_1x)+q(t,c_2x)$};
\node at (1.5,3.3) { $p^{(2)}(t,x)$};
\node at (2.5,2) { $p^{(2)}(t,x)$};
\node at (4.3,1.2) { $p^{(2)}(t,c_1x)$}; 
\node at (4.3,0.7) { $+q(t,c_2x)$}; 
\end{tikzpicture}

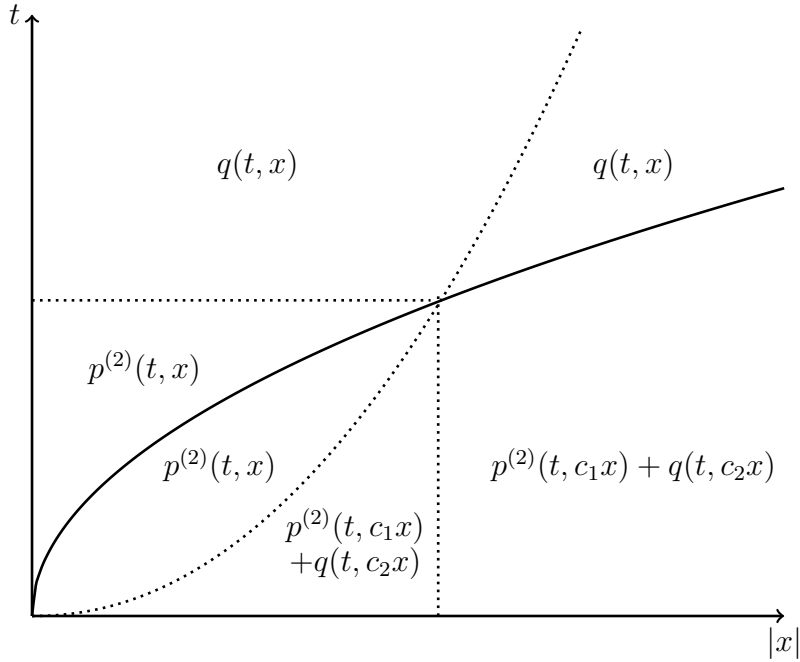
FIGURE1. Regions of heat kernel estimates for $p(t,x)$. Dotted line corresponds to $t=|x|^2$ and full line corresponds to $t\phi(|x|^{-2})=1$.\\
Note that heat kernel estimate is different to \cite[Theorem 1.4]{CK}, when $t\phi(|x|^{-2})\leq 1$ and $|x|>1$. $p^{(2)}(t,c_1x)$ additionally appear in our case.

\section{Dirichlet heat kernel estimates in $C^{1,1}$-open sets}
\noindent Recall that $p_D(t,x,y)$ is the transition density of $X^D$. In this section we obtain the sharp estimates of $p_D(t,x,y)$ in $C^{1,1}$-open sets.
\subsection{Lower bounds}  In this subsection we derive the lower bound estimate on $p_D(t,x,y)$ when $D$ is a $C^{1,1}$-open set. When $D$ is unbounded, we assume that the path distance in each connected component of $D$ is comparable to the Euclidean distance. Since the proofs are almost identical, we will provide a proof when $D$ is a bounded $C^{1,1}$-open set. We will use some relation between killed subordinate Brownian motions and subordinate killed Brownian motions. 

Let $T_t$ be a subordinator whose Laplace exponent $\phi$ is given by (\ref{lapex}). Then $t+T_t$ is a subordinator which has the same law as $S_t$. So $\{X_t;t\geq 0\}$ starting from $x$ has the same distribution as $\{B_{t+T_t};t \geq 0\}$ starting from $x$. Suppose that $U$ is an open subset of $\R^d$. We denote by $B^U$ the part process of $B$ killed upon leaving $U$. The process $\{Z_t^U;t\geq 0\}$ defined by $Z_t^U=B^U_{t+T_t}$ is called a subordinate killed Brownian motion in $U$. Let $q_U(t,x,y)$ be the transition density of $Z^U$. Denote by $\zeta^{Z,U}$ the lifetime of $Z^U$. Clearly, $Z_t^U =B_{t+T_t}$ for every $t\in [0, \zeta^{Z,U})$. Therefore we have
$$
p_U(t,z,w)\geq q_U(t,z,w) \quad \text{for}\quad (t,z,w) \in (0,\infty)\times U\times U.
$$

In the next proposition we will use \cite[Proposition 2.4]{M}.
Note that there is  a typo in \cite[Proposition 2.4]{M}. 
$\alpha\phi^{-1}(\beta^{-1})$  in the display there should be $\alpha\phi^{-1}(\beta t^{-1})$.

\begin{prop}\label{scces}
Suppose that $D$ is a $C^{1,1}$-open set in $\R^d$ with characteristics $(R_0,\Lambda_0)$. If $D$ is bounded, we assume that $\phi$ satisfies $L_{a}(\gamma, C_L)$ for some $a>0$. If $D$ is unbounded, we assume that $\phi$ satisfies $L_{0}(\gamma, C_L)$ and the path distance in each connected component of $D$ is comparable to the Euclidean distance with characteristic $\lambda_0$. For any $T>0$ there exist positive constants $c_1=c_1(R_0,\Lambda_0,\lambda_0,T, \phi)$ and $c_2=(R_0, \Lambda_0, \lambda_0)$ such that for all $t \in (0,T]$ and $x,y$ in the same connected component of $D$, 
$$
p_D(t,x,y) \geq c_1\bigg(1 \wedge \frac{\delta_D(x)}{\sqrt{t}}\bigg)\bigg(1 \wedge \frac{\delta_D(y)}{\sqrt{t}}\bigg) \phi^{-1}(t^{-1})^{d/2}e^{-c_2|x-y|^2\phi^{-1}(t^{-1})}.
$$
\end{prop}
\pf
Suppose that $x$ and $y$ are in the same component of $D$ and $\rho \in (0,1)$ is the constant in \cite[Proposition 2.4]{M} for $\alpha=2$ and $\beta=1$. Without loss of generality, let $T\geq 1$ and  $a$ satisfies $T<\rho\phi(a)^{-1}$ by Remark \ref{msc}. Let $\tilde{p}_D(t,z,w)$ be the transition density of $B^D$ (killed Brownian motion) and $q_D(t,x,y)$ be the transition density of $B^D_{S_t}$ (subordinate killed Brownian motion). By \cite[Theorem 3.3]{Ch} (see also \cite[Theorem 1.2]{Zh}) (where the comparability condition on the path distance in each component of $D$ with the Euclidean distance is used if $D$ is unbounded), there exists positive constants $c_3=c_3(R_0, \Lambda_0, \lambda_0,T, \phi)$ and $c_4=c_4(R_0,\Lambda_0,\lambda_0)$ such that for any $(s,z,w)\in (0,\phi^{-1}(\rho T^{-1})^{-1}] \times D \times D$,
$$
\tilde{p}_D(s,z,w) \geq c_3\Big(1\wedge \frac{\delta_D(z)}{\sqrt{s}}\Big)\Big(1\wedge \frac{\delta_D(w)}{\sqrt{s}}\Big)s^{-d/2}e^{-c_4|x-w|^2/2}.
$$
(Although not explicitly mentioned in \cite{Ch}, a careful examination of the proofs in \cite{Ch} reveals that the constants $c_3$ and $c_4$ in the above lower bound estimate can be chosen to depend only on $(R_0,\Lambda_0,\lambda_0,T, \phi)$ and $(R_0,\Lambda_0,\lambda_0)$, respectively.)\\
We have that for $0<t\leq T$,
\begin{align*}
p_D(t,x,y) &\geq q_D(t,x,y) \\
&=\int_{(0,\infty)}\tilde{p}_D(s,x,y)\P(S_t \in ds)\\
&\geq c_3\int_{[2^{-1}\phi^{-1}(t^{-1})^{-1},\phi^{-1}(\rho t^{-1})^{-1}]}\Big(1\wedge \frac{\delta_D(x)}{\sqrt{s}}\Big)\Big(1\wedge \frac{\delta_D(y)}{\sqrt{s}}\Big)s^{-d/2}e^{-c_4\frac{|x-y|^2}{s}}\P(S_t \in ds)\\
&\geq c_3\Bigg(1\wedge \frac{\delta_D(x)}{\sqrt{\phi^{-1}(\rho t^{-1})^{-1}}}\Bigg)\Bigg(1\wedge \frac{\delta_D(y)}{\sqrt{\phi^{-1}(\rho t^{-1})^{-1}}}\Bigg)\phi^{-1}(\rho t^{-1})^{d/2}e^{-2c_4|x-y|^2\phi^{-1}(t^{-1})}\\
&\quad \times \P(2^{-1}\phi^{-1}(t^{-1})^{-1}\leq S_t \leq \phi^{-1}(\rho t^{-1})^{-1}).
\end{align*}
Since $0<t<\rho\phi(a)^{-1}$, using the condition $L_a(\gamma,C_L)$ on $\phi$ (also see \cite[Lemma 3.1(ii)]{M}), we have
$$
\phi^{-1}(\rho t^{-1})=\phi^{-1}(t^{-1})\frac{\phi^{-1}(\rho t^{-1})}{\phi^{-1}(t^{-1})} \geq C_L^{1/\gamma}\rho^{1/\gamma}\phi^{-1}(t^{-1}).
$$
Using this and \eqref{phisc}, we have 
$$
\phi^{-1}(\rho t^{-1}) \geq C_L^{2/\gamma}\rho^{1/\gamma}T^{1-1/\gamma}\phi^{-1}(1)t^{-1}.
$$
Using the last two displays and \cite[Proposition 2.4]{M} we get
\begin{align*}
p_D(t,x,y) &\geq c_5\bigg(1 \wedge \frac{\delta_D(x)}{\sqrt{t}}\bigg)\bigg(1 \wedge \frac{\delta_D(y)}{\sqrt{t}}\bigg) \phi^{-1}(\rho t^{-1})^{d/2}e^{-2c_4|x-y|^2\phi^{-1}(t^{-1})}\\
&\quad \times \P(2^{-1}\phi^{-1}(t^{-1})^{-1}\leq S_t \leq \phi^{-1}(\rho t^{-1})^{-1})\\
&\geq c_6\bigg(1 \wedge \frac{\delta_D(x)}{\sqrt{t}}\bigg)\bigg(1 \wedge \frac{\delta_D(y)}{\sqrt{t}}\bigg) C_L^{d/2\gamma}\rho^{d/2\gamma} \phi^{-1}(t^{-1})^{d/2}e^{-2c_4|x-y|^2\phi^{-1}(t^{-1})}\tau,
\end{align*}
where $\tau$ is the constant in \cite[Proposition 2.4]{M} for $\alpha=2$ and $\beta=1$. \qed\\

 \begin{lemma}\label{L:5.1}
 For any positive constants $ a, b$ and $T$, there exists $c>0$ such that for all $z \in \R^d$ and $0<t\leq T$,
 $$
 \inf_{y \in B(z,at^{1/2}/2)} \P^y(\tau_{B(z,at^{1/2})}>bt)\geq c
 $$
 \end{lemma}
\pf
See \cite[Lemma 2.3]{CKS7}.
\qed

Although the proof of the following Lemma is similar to that of Lemma 2.4 of \cite{CKS7}, we give the proof again to make the paper self-contained.
\begin{lemma}\label{leb}
Suppose $H$ satisfies $L_a(\gamma, C_L)$ and $U_a(\delta, C_U)$ with $\delta <2$ for some $a>0$ ($L_0(\gamma,C_L)$ and $U_0(\delta,C_U)$, respectively). Then 
for every $T>0, M>0$ and $b>0$ 
there exists $c>0$ such that 
we have that for all $t \in(0,T]$, and $u,v \in \R^d$ satisfying $|u-v|\leq M/2$ ($u,v\in \R^d$, respectively)
$$
p_E(t,u,v) \geq c(t^{-d/2}\wedge t|u-v|^{-d}H(|u-v|^{-2}))
$$
where $E:=B(u, bt^{1/2}) \cup B(v, bt^{1/2})$.
\end{lemma}

\pf We fix $b>0$ and $u,v \in \R^d$ satisfying $|u-v|\leq M/2$, and let $r_t:=bt^{1/2}$. If $|u-v| \leq r_t/2$, by \cite[Lemma 2.1]{CKS7} (with $\sqrt{\lambda}=r_t$ and $D=B(0,1)$), 
\begin{align*}
p_E(t,u,v) &\geq \inf_{|z|<r_t/2}p_{B(0,r_t)}(t,0,z)=\inf_{|z|<r_t/2}p_{B(0,r_t)}(r_t^2(t/r_t^2),0,z)\\
&\geq c_1t^{-d/2}\bigg(1\wedge \frac{r_t}{\sqrt{t}}\bigg)\bigg(1\wedge \frac{r_t}{2\sqrt{t}}\bigg)e^{-c_2r_t^2/t}\geq c_3t^{-d/2}.
\end{align*}
If $|u-v|\geq r_t/2$, since the distance between $B(u,r_t/8)$ and $B(v,r_t/8)$ is at least $r_t/4$, we have by the strong Markov property and the L\'evy system of $X$ in (\ref{levy}) that
\begin{align*}
&p_E(t,u,v) \\
&\geq \E_u[p_E(t-\tau_{B(u,r_t/8)}, X_{\tau_{B(u,r_t/8)}},v): \tau_{B(u,r_t/8)}<t, X_{\tau_{B(u,r_t/8)}} \in B(v, r_t/8)] \\
& = \int_0^t \Big( \int_{B(u,r_t/8)}p_{B(u,r_t/8)}(s,u,w)\Big( \int_{B(v,r_t/8)}J(w,z)p_E(t-s,z,v)dz\Big) dw\Big)ds \\
&\geq \Big(\inf_{w\in B(u,r_t/8),z\in B(v,r_t/8)}J(w,z)\Big)\int_0^t \P_u(\tau_{B(u,r_t/8)}>s)\Big(\int_{B(v,r_t/8)}p_E(t-s,z,v)dz\Big)ds \\
&\geq \P_u(\tau_{B(u,r_t/8)}>t)\Big(\inf_{w\in B(u,r_t/8),z\in B(v,r_t/8)}J(w,z)\Big)\int_0^t \int_{B(v,r_t/8)}p_{B(v,r_t/8)}(t-s,z,v)dzds\\
&=\P_0(\tau_{B(0,r_t/8)}>t) \Big(\inf_{w\in B(u,r_t/8),z\in B(v,r_t/8)}j(|w-z|)\Big)\int_0^t \P_0 (\tau_{B(0,r_t/8)}>s) \\
&\geq t(\P_0(\tau_{B(0,r_t/8)}>t) )^2 \Big(\inf_{w\in B(u,r_t/8),z\in B(v,r_t/8)}j(|w-z|)\Big) \\
&\geq c_5t\Big(\inf_{w\in B(u,r_t/8),z\in B(v,r_t/8)}j(|w-z|)\Big).
\end{align*}
In the last inequality we have used Lemma \ref{L:5.1}. Note that if $w\in B(u,r_t/8)$ and $z\in B(v,r_t/8)$, then 
$$
|w-z| \leq |u-w|+|u-v|+|v-z| \leq |u-v| +\frac{r_t}{4} \leq 2|u-v|.
$$
Thus using (\ref{jl}) and Lemma \ref{lem:bf} we have
\begin{align*}
p_E(t,u,v)&\geq c_5tj(2|u-v|)\geq c_6t2^{-d}|u-v|^{-d}H(|u-v|^{-2}/4)\\
&\geq 2^{-d-4}c_6t|u-v|^{-d}H(|u-v|^{-2}).
\end{align*}
\qed

The next lemma say that if $x$ and $y$ are far away, the jumping kernel component dominates the Gaussian component and another off-diagonal estimate component.
\begin{lemma}\label{L:5.5}
Suppose $\phi$ satisfies $L_a(\gamma, C_L)$ for some $a\geq 0$. For any given positive constants $c_1,c_2,R$ and $T$, there is a positive constant $c_3=c_3(R,T,c_1,c_2)$ so that
$$
t^{-d/2}e^{-r^2/(c_1t)}+\phi^{-1}(t^{-1})^{d/2}e^{-c_2r^2\phi^{-1}(t^{-1})} \leq c_3tr^{-d}H(r^{-2})$$
for every $r \geq R$ and $t\in(0,T]$.
\end{lemma}
\pf
By Lemma \ref{lem:bf} there exist $c_4>0$ and $c_5>0$ such that 
$$
r^{-d}H(r^{-2}) \geq r^{-d-4} \geq c_4e^{-c_5r} \quad \text{for every} \quad r>1.
$$
For $r>1 \vee (2c_1c_5T)\vee \frac{2c_5}{c_2}\phi^{-1}(T^{-1})^{-1}$ and $t \in (0,T]$, we have following inequalities  
$$r^2/(2c_1t)>c_5r, \quad c_2r^2\phi^{-1}(t^{-1})/2 >c_5r,$$
\begin{align*}
t^{-d/2-1}e^{-r^2/(2c_1t)} &\leq t^{-d/2-1}e^{-1/(2c_1t)} \leq \sup_{0<s\leq T}s^{-d/2-1}e^{-1/(2c_1s)} =:c_6 <\infty ,
\end{align*}
and 
\begin{align*}
\phi^{-1}(t^{-1})^{d/2}t^{-1}&e^{-c_2r^2\phi^{-1}(t^{-1})/2} \\
\leq& \sup_{0<s\leq T} \phi^{-1}(s^{-1})^{d/2}s^{-1}e^{-c_2\phi^{-1}(s^{-1})/2} \\
\leq &T^{1/\gamma -1}C_L^{-1/\gamma}\phi^{-1}(1)^{-1}\sup_{0<s\leq T}\phi^{-1}(s^{-1})^{d/2+1}e^{-c_2\phi^{-1}(s^{-1})/2}
=:c_7 <\infty.
\end{align*}
In the last inequality we have used \eqref{phisc}. (Without loss of generality, we can assume that $a\leq \phi^{-1}(T^{-1})$.)
Therefore when $r>1 \vee (2c_1c_5T) \vee \frac{2c_5}{c_2}\phi^{-1}(T^{-1})^{-1}$ and $t \in (0,T]$, we have
$$
t^{-d/2}e^{-r^2/(c_1t)} \leq c_6te^{-r^2/(2c_1t)} \leq c_6te^{-c_5r} \leq (c_6/c_4)tr^{-d}H(r^{-2})
$$
and 
$$
\phi^{-1}(t^{-1})^{d/2}e^{-c_2r^2\phi^{-1}(t^{-1})} \leq c_7te^{-c_2r^2\phi^{-1}(t^{-1})/2} \leq c_7te^{-c_5r} \leq  (c_7/c_4)tr^{-d}H(r^{-2}).
$$
When $R\leq r \leq 1 \vee (2c_1c_5T) \vee \frac{2c_5}{c_2}\phi^{-1}(T^{-1})^{-1}$ and $t\in (0,T]$, clearly
$$
t^{-d/2}e^{-r^2/(c_1t)} \leq t\Bigg(\sup_{s\leq T}s^{-d/2-1}e^{-R^2/(c_1s)}\Bigg) \leq c_8tr^{-d}H(r^{-2}) \quad
$$
and 
$$
\phi^{-1}(t^{-1})^{d/2}e^{-c_2r^2\phi^{-1}(t^{-1})}  \leq t\Bigg(\sup_{s\leq T}\phi^{-1}(s^{-1})^{d/2}s^{-1}e^{-c_2R^2\phi^{-1}(s^{-1})}\Bigg)\leq c_9tr^{-d}H(r^{-2}).
$$
\qed\\
{\bf Proof of Theorem \ref{T:3} (2) and Theorem \ref{T:4} (2)}. Since two proofs are almost identical, we just prove Theorem \ref{T:3} (2).
First note that the distance between two distinct connected components of $D$ is at least $R_0$. 
Since $D$ is a $C^{1,1}$ open set, it satisfies the uniform interior ball condition with radius $r_0=r_0(R_0, \Lambda_0)\in (0,R_0]$: there exists $r_0=r_0(R_0, \Lambda_0)\in (0,R_0]$ such that for any $x\in D$ with $\delta_D(x)<r_0$, there are $z_x\in \partial D$ so that $|x-z_x|=\delta_D(x)$ and that $B(x_0,r_0)\subset D$ for $x_0=z_x+r_0(x-z_x)/|x-z_x|$. Set $T_0=(r_0/4)^2$. Using such uniform interior ball condition, 
by considering the cases $\delta_D(x)<r_0$ and  $\delta_D(x)>r_0$, there exists $L=L(r_0)>1$ such that, for all $t\in (0,T_0]$ and $x,y \in D$, we can choose $\xi^t_x \in D\cap B(x,L\sqrt{t})$ and $\xi^t_y \in D \cap B(y,L\sqrt{t})$ so that $B(\xi^t_x,2\sqrt{t})$ and $B(\xi^t_y,2\sqrt{t})$ are subsets of the connected components of $D$ that contains $x$ and $y$, respectively. 

We first consider the case $t\in(0,T_0]$.
Note that by the semigroup property,
\begin{equation}\label{sp}
p_D(t,x,y)\geq \int_{B(\xi^t_x,\sqrt{t})}\int_{B(\xi^t_y,\sqrt{t})} p_D(t/3,x,u)p_D(t/3,u,v)p_D(t/3,v,y)dudv.
\end{equation}
For $u\in B(\xi^t_x,\sqrt{t})$, we have 
$$\delta_D(u) \geq \sqrt{t} \quad \text{and} \quad |x-u|\leq |x-\xi^t_x|+|\xi^t_x-u| \leq L\sqrt{t}+\sqrt{t}=(L+1)\sqrt{t}.
$$
Thus by \cite[Lemma 2.1]{CKS7}, for $t\in(0,T_0]$,
\begin{eqnarray}
\int_{B(\xi^t_x, \sqrt{t})}p_D(t/3,x,u)du \geq c_3\bigg(1\wedge \frac{\delta_D(x)}{\sqrt{t}}\bigg)\int_{B(\xi^t_x,\sqrt{t})}\bigg(1\wedge \frac{\delta_D(u)}{\sqrt{t}}\bigg)t^{-d/2}e^{-c_4|x-u|^2/t}du\nonumber\\
\geq c_3\bigg(1\wedge \frac{\delta_D(x)}{\sqrt{t}}\bigg)t^{-d/2}e^{-c_4(L+1)^2}|B(\xi^t_x,\sqrt{t})| \geq c_5\bigg(1\wedge \frac{\delta_D(x)}{\sqrt{t}}\bigg). \label{e:x}
\end{eqnarray}
Similarly, for $t\in(0,T_0]$ and $v\in B(\xi^t_y,\sqrt{t})$,
\begin{equation}\label{e:y}
\int_{B(\xi^t_y, \sqrt{t})}p_D(t/3,y,v)dv \geq c_5\bigg(1\wedge \frac{\delta_D(y)}{\sqrt{t}}\bigg).
\end{equation}
 Using (\ref{sp}), Proposition \ref{leb}, symmetry and (\ref{e:x})--(\ref{e:y}), we have
\begin{align}\label{pDlower1}
&p_D(t,x,y)\nn\\
&\geq \int_{B(\xi^t_y,\sqrt{t})}\int_{B(\xi^t_x,\sqrt{t})}p_D(t/3,x,u)p_{B(u,\sqrt{t}/2)\cup B(v,\sqrt{t}/2)}(t/3,u,v)p_D(t/3,v,y)dudv\nn\\
&\geq c_6\int_{B(\xi^t_y,\sqrt{t})}\int_{B(\xi^t_x,\sqrt{t})}p_D(t/3,x,u)\bigg(t^{-d/2}\wedge \big(\frac{tH(|u-v|^{-2})}{|u-v|^{d}}\big)\bigg)p_D(t/3,v,y)dudv\nn\\
&\geq c_6\Bigg(\inf_{(u,v)\in B(\xi^t_x,\sqrt{t})\times B(\xi^t_y,\sqrt{t})}\bigg(t^{-d/2}\wedge \big(\frac{tH(|u-v|^{-2})}{|u-v|^{d}}\big)\bigg)\Bigg)\nn\\
&\quad \times \int_{B(\xi^t_y,\sqrt{t})}\int_{B(\xi^t_x,\sqrt{t})}p_D(t/3,x,u)p_D(t/3,v,y)dudv\nn\\
&\geq c_6c_5^2\Bigg(\inf_{(u,v)\in B(\xi^t_x,\sqrt{t}) \times B(\xi^t_y,\sqrt{t})}\bigg(t^{-d/2}\wedge \big(\frac{tH(|u-v|^{-2})}{|u-v|^{d}}\big)\bigg)\Bigg) \bigg(1 \wedge \frac{\delta_D(x)}{\sqrt{t}}\bigg)\bigg(1 \wedge \frac{\delta_D(y)}{\sqrt{t}}\bigg).
\end{align}
Suppose that $|x-y| \geq \sqrt{t}/8$ and $t \in(0,T_0]$. Then  we have that for $(u,v) \in B(\xi^t_x,\sqrt{t}) \times B(\xi^t_y,\sqrt{t})$,
\begin{align*}
&|u-v|\leq |u-\xi^t_x|+|\xi^t_x-x|+|x-y|+|y-\xi^t_y|+|\xi^t_y-v|\\
&\leq 2(1+L)\sqrt{t}+|x-y|\leq (16(1+L)|x-y|).
\end{align*}
Thus using Lemma \ref{lem:bf} we have 
\begin{align*}
\inf_{(u,v)\in B(\xi^t_x,\sqrt{t}) \times B(\xi^t_y,\sqrt{t})}\bigg(t^{-d/2}\wedge \big(\frac{tH(|u-v|^{-2})}{|u-v|^{d}}\big)\bigg)
\geq c_7\bigg(t^{-d/2}\wedge \big(\frac{tH(|x-y|^{-2})}{|x-y|^{d}}\big)\bigg).
\end{align*}
Therefore,  for $|x-y| \geq \sqrt{t}/8$ and $t \in(0,T_0]$
\begin{equation}\label{e:l1}
p_D(t,x,y)\geq c_8\bigg(1 \wedge \frac{\delta_D(x)}{\sqrt{t}}\bigg)\bigg(1 \wedge \frac{\delta_D(y)}{\sqrt{t}}\bigg) \bigg(t^{-d/2}\wedge \big(\frac{tH(|x-y|^{-2})}{|x-y|^{d}}\big)\bigg).
\end{equation}

Using the inequality (\ref{e:l1}), we will obtain the sharp lower bound estimates by
considering  the following three cases. \\
{\it Case (1):} Suppose that $|x-y| \geq \sqrt{t}/8$, $t \in(0,T_0]$, and $x$ and $y$ are contained in same connected component of $D$. Combining with (\ref{e:l1}), Proposition \ref{scces}, and \cite[Lemma 2.1]{CKS7}, we conclude that
\begin{align}\label{e:l2}
p_D(t,x,y) &\geq c_9\bigg(1 \wedge \frac{\delta_D(x)}{\sqrt{t}}\bigg)\bigg(1 \wedge \frac{\delta_D(y)}{\sqrt{t}}\bigg)\nn\\
&\quad \times \left(\left(t^{-d/2}\wedge \big(\frac{tH(|x-y|^{-2})}{|x-y|^{d}}\big)\right)+\phi^{-1}(t^{-1})^{d/2}e^{-c_{10}|x-y|^2\phi^{-1}(t^{-1})}+t^{-d/2}e^{-\frac{|x-y|^2}{c_{11}t}}\right)\nonumber\\
& \geq c_{9}\bigg(1 \wedge \frac{\delta_D(x)}{\sqrt{t}}\bigg)\bigg(1 \wedge \frac{\delta_D(y)}{\sqrt{t}}\bigg) \nonumber\\
& \quad \times \left(t^{-d/2}\wedge (t^{-d/2}e^{-\frac{|x-y|^2}{c_{11}t}}+\frac{tH(|x-y|^{-2})}{|x-y|^{d}}+\phi^{-1}(t^{-1})^{d/2}e^{-c_{10}|x-y|^2\phi^{-1}(t^{-1})})\right).\nonumber\\
\end{align}
{\it Case (2):} Suppose that $|x-y| \geq \sqrt{t}/8$, $t \in(0,T_0]$, and $x$ and $y$ are contained in two distinct connected components of $D$. By (\ref{e:l1}) and Lemma \ref{L:5.5}, we have the same conclusion in (\ref{e:l2}).\\
{\it Case (3):} Suppose that $|x-y|<\sqrt{t}/8$ and $t\in (0,T_0]$. In this case $x$ and $y$ are in the same connected component. For $(u,v) \in B(\xi^t_x,\sqrt{t}) \times B(\xi^t_y,\sqrt{t})$,
$$
|u-v|\leq 2(1+L)\sqrt{t}+|x-y|\leq (2(1+L)+8^{-1})\sqrt{t}.
$$
Thus by \cite[Lemma 2.1]{CKS7}, we have that for every $(u,v) \in B(\xi^t_x,\sqrt{t}) \times B(\xi^t_y,\sqrt{t})$,
$$
p_D(t/3,u,v) \geq c_{12}\bigg(1\wedge \frac{\delta_D(u)}{\sqrt{t}}\bigg)\bigg(1\wedge\frac{\delta_D(v)}{\sqrt{t}}\bigg)t^{-d/2}e^{-c_{13}|u-v|^2/t}\geq c_{14}t^{-d/2}.
$$
Therefore by (\ref{sp})--(\ref{e:y}), for $t\leq T_0$,
{\setlength\arraycolsep{-3pt}
\begin{eqnarray}\label{e:l3}
&&p_D(t,x,y)\geq c_{14}c_5^2\bigg(1 \wedge \frac{\delta_D(x)}{\sqrt{t}}\bigg)\bigg(1 \wedge \frac{\delta_D(y)}{\sqrt{t}}\bigg) t^{-d/2}\nonumber\\
&&\geq c_{14}c_5^2\bigg(1 \wedge \frac{\delta_D(x)}{\sqrt{t}}\bigg)\bigg(1 \wedge \frac{\delta_D(y)}{\sqrt{t}}\bigg) \nonumber\\
&&\times \left(t^{-d/2} \wedge \left(t^{-d/2}e^{-\frac{|x-y|^2}{c_{15}t}}+\frac{tH(|x-y|^{-2})}{|x-y|^{d}}+\phi^{-1}(t^{-1})^{d/2}e^{-c_{16}|x-y|^2\phi^{-1}(t^{-1})}\right)\right).\nonumber\\
\end{eqnarray}}
Combining the above three cases, we get (\ref{dhkel}) for $t\in(0,T_0]$. When $T>T_0$ and $t \in(T_0,T]$, observe that $T_0/3 \leq t-2T_0/3 \leq T-2T_0/3 \leq (T/T_0-2/3)T_0$, that is, $t-2T_0/3$ is comparable to $T_0/3$ with some universal constants that depend only on $T$ and $T_0$. Using the inequality
\begin{align*}
&p_D(t,x,y)\geq \int_{B(\xi^{T_0}_x,\sqrt{T_0})}\int_{B(\xi^{T_0}_y,\sqrt{T_0})} p_D(T_0/3,x,u)p_D(t-2T_0/3,u,v)p_D(T_0/3,v,y)dudv\\
&\geq  \int_{B(\xi^{T_0}_x,\sqrt{T_0})}\int_{B(\xi^{T_0}_y,\sqrt{T_0})} p_D(T_0/3,x,u)p_{B(u,\sqrt{T_0}/2)\cup B(v,\sqrt{T_0}/2)}(t-2T_0/3,u,v)p_D(T_0/3,v,y)dudv
\end{align*}
instead of (\ref{sp}) and following the argument in (\ref{pDlower1}) and (\ref{e:l1}) we have
$$
p_D(t,x,y)\geq c_{17}\bigg(1 \wedge \frac{\delta_D(x)}{\sqrt{T_0}}\bigg)\bigg(1 \wedge \frac{\delta_D(y)}{\sqrt{T_0}}\bigg) \bigg(T_0^{-d/2}\wedge \big(T_0\frac{H(|x-y|^{-2})}{|x-y|^{d}}\big)\bigg).
$$ Consider the cases $|x-y|\geq \sqrt{T_0}/8$ and $|x-y|<\sqrt{T_0}/8$ separately and follow the above three cases. Then since $\frac{tT_0}{T}\leq T_0< t$ for $t\in[T_0,T]$, we can obtain (\ref{dhkel}) for $t\in [T_0,T]$ and hence for $t\in(0,T]$.
\qed
\subsection{Upper bounds} In this subsection we derive the upper bound estimate on $p_D(t,x,y)$ when $D$ is a $C^{1,1}$-open set(not necessarily bounded). We use the following lemma in \cite[Lemma 3.1]{CKS7}.
\begin{lemma}[{\cite[Lemma 3.1]{CKS7}}]\label{L:5.7}
Suppose that $U_1, U_3, E$ are open subsets of $\R^d$ with $U_1,U_3 \subset E$ and $\dist (U_1,U_3)>0$. Let $U_2:=E\setminus (U_1 \cup U_3)$. If $x\in U_1$ and $y\in U_3$, then for every $t>0$,
{\setlength\arraycolsep{-8pt}
\begin{eqnarray}
&&p_E(t,x,y) \leq \P_x\big(X_{\tau_{U_1}}\in U_2\big)\bigg(\sup_{s<t, z\in U_2}p_E(s,z,y)\bigg) \nonumber\\
&&\quad\quad\quad\quad\quad+ \int_0^t \P_x(\tau_{U_1}>s)\P_y(\tau_E > t-s)ds\bigg(\sup_{u\in U_1, z\in U_3}J(u,z)\bigg)\label{3.1} \\
&&\quad\quad\quad\quad\leq \P_x\big(X_{\tau_{U_1}}\in U_2\big)\bigg(\sup_{s<t, z\in U_2}p(s,z,y)\bigg) +(t\wedge \E_x[\tau_{U_1}])\bigg(\sup_{u\in U_1, z\in U_3}J(u,z)\bigg).\label{3.2}
\end{eqnarray}}
\end{lemma}
Note that by Remark \ref{R:L2.3} or Theorem \ref{T:2}, there exist positive constants $c, a_U$, and $C_2$ such that
\begin{equation}\label{uu}
p(t,x,y)\leq c\Big(t^{-d/2}\wedge \big(t^{-d/2}e^{-\frac{|x-y|^2}{C_2t}} + \frac{tH(|x-y|^{-2})}{|x-y|^d}+\phi^{-1}(t^{-1})^{d/2}e^{-a_U|x-y|^2\phi^{-1}(t^{-1})}\big)\Big).
\end{equation}
The boundary Harnack principle for subordinate Brownian motions with Gaussian components was proved 
in \cite{KSV1} for any $C^{1,1}$-open set, see \cite[Theorem 1.2]{KSV1}. 
In  \cite[Theorem 1.2]{KSV1}, it is assumed that $\phi$ is a complete Bernstein function and that the L\'evy density $\mu$ of $S$ satisfies growth condition near zero, i.e., for any $K>0$, there exists $c=c(K)>1$ such that $\mu(r)\leq c\mu(2r)$.

Note that in the proof of \cite[Theorem 1.2]{KSV1}, as a consequence of the growth condition of L\'evy density of $S$ and assumption that $\phi$ is a complete Bernstein function, in fact, the following conditions of L\'evy density $j$ of $X$ are actually used (see \cite[(2.7), (2.8)]{KSV1}):\\
for any $K>0$, there exists $c_1=c_1(K)>1$ such that
\begin{equation}\label{e:2.6}
j(r)\leq c_1j(2r), \quad \text{for}\quad r\in(0,K),
\end{equation}
and there exists $c_2>1$ such that 
\begin{equation}\label{e:2.7}
j(r)\leq c_2j(r+1), \quad \text{for}\quad r>1.
\end{equation}
If, instead of the assumption that $\phi$ is a complete Bernstein function and L\'evy density $\mu$ satisfies growth condition near zero,  we assume that $H$ satisfies $L_0(\gamma, C_L)$ and $U_0(\delta, C_U)$ with $\delta<2$, then by (\ref{jl}), (\ref{e:2.6}) and (\ref{e:2.7}) hold. Thus the boundary Harnack principle still hold. But if we assume that $H$ satisfies $L_a(\gamma, C_L)$ and $U_a(\delta, C_U)$ with $\delta <2$  for some $a>0$, then (\ref{e:2.7}) may not holds. Nonetheless, if we only consider harmonic functions not only vanishing continuously on $D^c \cap B(Q,r)$, $Q\in \partial D$, but also zero on $D^c$, then we don't need the condition (\ref{e:2.7}). Thus we have the following modified theorem. 
\begin{thm}\label{mbhp} Let $D$ is a $C^{1,1}$-open set in $\R^d$ with characteristics $(R_0, \Lambda_0)$. If $D$ is bounded, then we assume that $H$ satisfies $L_a(\gamma, C_L)$ and $U_a(\delta, C_U)$ with $\delta <2$ for some $a>0$. If $D$ is unbounded, then we assume that $H$ satisfies $L_0(\gamma, C_L)$ and $U_0(\delta, C_U)$ with $\delta <2$. Then there exists a positive constant $c=c(d, \Lambda_0, R_0)$ such that for $r\in(0,R_0]$, $Q\in \partial D$ and any nonnegative function $f$ in $\R^d$ which is harmonic in $D\cap B(Q,r)$ with respect to $X$, zero on $D^c$ and vanishes continuously on $D^c \cap B(Q,r)$, we have
\begin{equation}
\frac{f(x)}{\delta_D(x)}\leq c\frac{f(y)}{\delta_D(y)}\quad \text{for every}\quad x,y \in D\cap B(Q,r/2).
\end{equation}
\end{thm}
\pf Since we have explained before the statement of the theorem why theorem holds for the unbounded case,  
we will just prove the theorem  when $D$ is bounded and $H$ satisfies $L_a(\gamma, C_L)$ and $U_a(\delta, C_U)$ with $\delta <2$ for some $a>0$.

Let $\R^d_+=\{x=(x_1,...,x_{d-1},x_d):=(\tilde{x},x_d)\in \R^d:x_d>0\}$, $V$ is the potential measure of the ladder height process of $X_t^d$, where $X_t^d$ is d-th component of $X_t$, and $w(x)=V((x_d)^+)$. 
We first  show that \cite[Proposition 3.3]{KSV1} holds under our assumptions, i.e., we claim that for any positive constants $r_0$ and $M$, we have 
\begin{align}\label{prop3.3}
\sup_{x\in \R^d:0<x_d<M}\int_{B(x,r_0)^c\cap \R^d_+}w(y)j(|x-y|)dy<\infty,
\end{align}

Once we have \eqref{prop3.3}, 
then for $f$ which is harmonic in $D\cap B(Q,r)$ with respect to $X$, zero on $D^c$ and vanishes continuously on $D^c \cap B(Q,r)$ where $r\in (0,R_0]$ and $Q\in \partial D$, we can follow the proofs 
of \cite[Theorem 5.3  and Theorem 1.2]{KSV1} line by line without using \eqref{e:2.7}.

By \cite[Theorem 5, page 79]{B} and \cite[Lemma 2.1]{KSV2}, $V$ is absolutely continuous and has a continuous and strictly positive density $v$ such that $v(0+)=1$. It is also well known that $V$ is subadditive, i.e., $V(s+t)\leq V(s)+V(t)$, $s,t\in \R$ (See \cite[page 74]{B}.) and $V(\infty)=\infty$. Without loss of generality we assume that $\tilde{x}=0$. Note that for $0<x_d<M$ and $y\in B(x,r_0)^c$, 
\begin{align}\label{subadd}
w(y)=V((y_d)^+)\leq V(|y|) \leq V(M+|x-y|)\leq V(M)+V(|x-y|).
\end{align}
Let $L(r)=\int _{r}^{\infty}r^{d-1}j(r)dr$, then by \cite[(2.23)]{BGR}, $L(r)\leq c_1/V(r)^2$. Using (\ref{subadd}), the integration by parts and \cite[(2.23)]{BGR} twice, we have
\begin{align*}
&\sup_{x\in R^d:0<x_d<M}\int_{B(x,r_0)^c\cap \R^d_+}w(y)j(|x-y|)dy \\
&\leq\sup_{x\in R^d:0<x_d<M} \int_{B(x,r_0)^c}(V(M)+V(|x-y|)j(|x-y|)dy\\
&\leq c_2\int_{r_0}^{\infty}(V(M)+V(r))r^{d-1}j(r)dr =c_2L(r_0)(V(M)+V(r_0))+c_2\int _{r_0}^{\infty}V'(r)L(r)dr\\
&\leq c_2L(r_0)(V(M)+V(r_0))+c_3\int _{r_0}^{\infty}\frac{V'(r)}{V(r)^2}dr\\
&=c_2L(r_0)(V(M)+V(r_0))+\frac{c_3}{V(r_0)} <\infty.
\end{align*}
We have proved \eqref{prop3.3}.
\qed

For the remainder of this section, we follow  proofs of \cite[Proposition 3.2 and Theorem 1.1(i)]{CKS7}.
First note that for $C^{1,1}$-open set $D\subset \R^d$ with characteristics $(R_0, \Lambda_0)$, there exists $r_0=r_0(R_0,\Lambda_0) \in (0,R_0]$ such that $D$ satisfies the uniform interior and uniform exterior ball conditions with radius $r_0$. We will use such $r_0>0$ in the proof of the next proposition and Theorem \ref{T:3} (1).
\begin{prop}\label{P:3.2}
Let $D$ is a $C^{1,1}$-open set in $\R^d$ with characteristics $(R_0, \Lambda_0)$. If $D$ is bounded, then we assume that $H$ satisfies $L_a(\gamma, C_L)$ and $U_a(\delta, C_U)$ with $\delta <2$ for some $a>0$. If $D$ is unbounded, then we assume that $H$ satisfies $L_0(\gamma, C_L)$ and $U_0(\delta, C_U)$ with $\delta <2$. For every $T>0$, there exists $c>0$ such that for all $(t,x,y)\in(0,T] \times D \times D$,
\begin{align}\label{e:upro}
&p_D(t,x,y) \leq c\bigg(1 \wedge \frac{\delta_D(x)}{\sqrt{t}}\bigg)\nn\\
&\qquad\times \Big(t^{-d/2}\wedge \big(t^{-d/2}e^{-\frac{|x-y|^2}{4C_2t}} + \frac{tH(|x-y|^{-2})}{|x-y|^{d}}+\phi^{-1}(t^{-1})^{d/2}e^{-\frac{a_U}{4}|x-y|^2\phi^{-1}(t^{-1})}\big)\Big),
\end{align}
where the constants $C_2$ and $a_U$ are from (\ref{uu}).
\end{prop}

\pf
We prove the proposition for the case  that $D$ is bounded, $H$ satisfies $L_a(\gamma, C_L)$ and $U_a(\delta, C_U)$ with $\delta <2$ for some $a>0$ only, because  the proof of the other case is almost identical.

 Fix $T>0$ and $t\in (0,T]$. Let $x,y \in D$. We just consider the case $\delta_D(x)<r_0\sqrt{t}/(16\sqrt{T}) \leq r_0/16$, if not, we can directly obtain \eqref{e:upro} by \eqref{uu}. Choose $x_0 \in \partial D$ and $x_1\in D$ such that $\delta_D(x)=|x-x_0|$ and 
$
x_1=x_0+\frac{r_0\sqrt{t}}{16\sqrt{T}}{\bf n}(x_0),
$ 
respectively,
where ${\bf n}(x_0)=(x-x_0)/|x_0-x|$ be the unit inward normal of $D$ at the boundary point $x_0$.
Define
$$
U_1:=B(x_0,r_0\sqrt{t}/(8\sqrt{T}))\cap D.
$$
Since (\ref{prop3.3}) holds, by \cite[Lemma 4.3]{KSV1}
\begin{equation}\label{3.6}
\E_x[\tau_{U_1}]\leq c_1\sqrt{t}\delta_D(x).
\end{equation}
Using Theorem \ref{mbhp} and $\delta_D(x_1)=\frac{r_0\sqrt{t}}{16\sqrt{T}}$, we have
\begin{align}\label{3.5}
\P_x(X_{\tau_{U_1}} \in D \setminus U_1) \leq c_2\P_{x_1}(X_{\tau_{U_1}} \in D \setminus U_1)\frac{\delta_D(x)}{\delta_D(x_1)} \leq c_2\frac{16\sqrt{T}\delta_D(x)}{r_0\sqrt{t}}\leq c_3\bigg(1\wedge \frac{\delta_D(x)}{\sqrt{t}}\bigg).
\end{align}
Thus by (\ref{3.5}) and (\ref{3.6}) we have
\begin{align}\label{3.7}
\P_x(\tau_D >t/2) & \leq \P_x(\tau_{U_1} >t/2) + \P_x\big(X_{\tau_{U_1}}\in D\setminus U_1\big) \nonumber\\
&\leq \bigg(1\wedge \Big(\frac{2}{t}\E_x[\tau_{U_1}]\Big)\bigg)+\P_x(X_{\tau_{U_1}} \in D \setminus U_1) \leq c_4\bigg(1 \wedge \frac{\delta_D(x)}{\sqrt{t}}\bigg).
\end{align}

Now we estimate $p_D(t,x,y)$ considering two cases separately. 
Let $c_5:= (d/2)\vee $ $((dC_L^{-1/\gamma}$ $T^{1/\gamma-1}\phi^{-1}(1)^{-1})/(2C_2a_U)) \vee (r^2_0/(4C_2T))$ where  $a_U$ and $C_2$ are the constants in \eqref{uu}.

\noindent
{\it Case (1):} $|x-y| \leq 
2(C_2c_5)^{1/2} \sqrt{t}$. By the semigroup property, 
$$
p_D(t,x,y) =\int_D p_D(t/2,x,z)p_D(t/2,z,y)dz \leq \Bigg(\sup_{z,w \in D}p(t/2,z,w)\Bigg)\int_D p_D(t/2,x,z)dz.
$$
Using Theorem \ref{T:2} and (\ref{3.7}) in the above display, we have
\begin{align*}
p_D(t,x,y) \leq c_6(t/2)^{-d/2}\P_x(\tau_D >t/2)\leq c_6c_42^{d/2}t^{-d/2}\bigg(1\wedge \frac{\delta_D(x)}{\sqrt{t}}\bigg).
\end{align*}
Since
$
|x-y|^2/(4C_2t) \leq c_5,
$
we have
\begin{align}\label{ndu1}
p_D(t,x,y) \leq c_6c_42^{d/2}e^{c_5}t^{-d/2}e^{-|x-y|^2/(4C_2t)}\bigg(1\wedge \frac{\delta_D(x)}{\sqrt{t}}\bigg).
\end{align}

\noindent
{\it Case (2):} $|x-y| \geq 2(C_2c_5)^{1/2} \sqrt{t}$. Define
\begin{equation}\label{3.9}
U_3:=\{z\in D : |z-x|>|x-y|/2\} \quad \text{and} \quad U_2:=D\setminus(U_1 \cup U_3).
\end{equation}
For $z\in U_2$,
\begin{equation}\label{3.12}
\frac{3}{2}|x-y|\geq |x-y|+|x-z|\geq |z-y|\geq |x-y|-|x-z|\geq \frac{|x-y|}{2}.
\end{equation}

By our choice of $c_5$ and $\eqref{phisc}$ (we can assume that $a\leq \phi^{-1}(T^{-1})$ by Remark \ref{msc}), we have $t\leq |x-y|^2/(2dC_2)$ and $\phi^{-1}(t^{-1})^{-1}\leq C_L^{-1/\gamma}T^{1/\gamma-1}t\phi^{-1}(1)^{-1} \leq \frac{a_U|x-y|^2}{2d}$. Using this and the fact that $s \to s^{-d/2}e^{-\beta/s}$ is increasing on the interval $(0,2\beta/d]$, we get that for $s \le t$
$$s^{-d/2}e^{-\frac{|x-y|^2}{4C_2s}}+ \phi^{-1}(s^{-1})^{d/2}e^{-\frac{a_U}{4}|x-y|^2\phi^{-1}(s^{-1})}  \le 
 t^{-d/2}e^{-\frac{|x-y|^2}{4C_2t}}+  \phi^{-1}(t^{-1})^{d/2}e^{-\frac{a_U}{4}|x-y|^2\phi^{-1}(t^{-1})}.$$

Thus by (\ref{uu}) and \eqref{3.12},
\begin{align}
&\sup_{s\leq t, z \in U_2}p(s,z,y)\nonumber \\
&\leq c_7\sup_{s\leq t, |z-y|\geq|x-y|/2}\left( s^{-d/2}e^{-\frac{|z-y|^2}{C_2s}} + s\frac{H(|z-y|^{-2})}{|z-y|^{d}}+\phi^{-1}(s^{-1})^{d/2}e^{-a_U|z-y|^2\phi^{-1}(s^{-1})}\right)\nonumber\\
&\leq c_7\sup_{s\leq t}\left(s^{-d/2}e^{-\frac{|x-y|^2}{4C_2s}} + 2^ds\frac{H(4|x-y|^{-2})}{|x-y|^d}+\phi^{-1}(s^{-1})^{d/2}e^{-\frac{a_U}{4}|x-y|^2\phi^{-1}(s^{-1})}\right)\nonumber\\
&\leq c_7\left(t^{-d/2}e^{-\frac{|x-y|^2}{4C_2t}}  + 2^{d+4}t\frac{H(|x-y|^{-2})}{|x-y|^{d}}+\phi^{-1}(t^{-1})^{d/2}e^{-\frac{a_U}{4}|x-y|^2\phi^{-1}(t^{-1})}\right)\nonumber\\
&\leq c_8\left(\left(t^{-d/2}\right)\wedge \left(t^{-d/2}e^{-\frac{|x-y|^2}{4C_2t}}  +t\frac{H(|x-y|^{-2})}{|x-y|^{d}}+\phi^{-1}(t^{-1})^{d/2}e^{-\frac{a_U}{4}|x-y|^2\phi^{-1}(t^{-1})}\right)\right).
\label{3.13}
\end{align}
For the last inequality, we argue as follows: by the proof of \cite[Corollary 1.3]{M}
\begin{align}\label{e:weakupper}
2^{d+4}t\frac{H(|x-y|^{-2})}{|x-y|^{d}}+\phi^{-1}(t^{-1})^{d/2}e^{-\frac{a_U}{4}|x-y|^2\phi^{-1}(t^{-1})}  \le c_{9}  t\frac{\phi(|x-y|^{-2})}{|x-y|^{d}}. 
\end{align}
On the other hand, by Lemma \ref{lem:bf}
$$ 
  t\frac{\phi(|x-y|^{-2})}{|x-y|^{d}} \le   t\frac{\phi( (4C_2c_5)^{-1} t^{-1})}{|x-y|^{d}} 
 \leq \frac{T\phi((4C_2c_5)^{-1}T^{-1}) }{|x-y|^{d}} \leq c_{10}t^{-d/2}.$$
 Therefore 
\begin{align}\label{3.14}
2^{d+4}t\frac{H(|x-y|^{-2})}{|x-y|^{d}}+\phi^{-1}(t^{-1})^{d/2}e^{-\frac{a_U}{4}|x-y|^2\phi^{-1}(t^{-1})} \le c_{11} t^{-d/2}.
\end{align}

For $u\in U_1$ and $z\in U_3$, since $|z-x|/2>|x-y|/4\geq \frac{r_0\sqrt{t}}{4\sqrt{T}}$,
$$
|u-z|\geq |z-x|-|x-x_0|-|x_0-u|\geq |z-x|-\frac{r_0\sqrt{t}}{4\sqrt{T}}\geq \frac{|z-x|}{2} >\frac{|x-y|}{4} \geq \frac{r_0\sqrt{t}}{4\sqrt{T}}.
$$
Thus $\dist(U_1,U_3)>0$ and by (\ref{jl}) and Lemma \ref{lem:bf}, we have
\begin{align}
\sup_{u\in U_1, z\in U_3}J(u,z)\leq \sup_{|u-z|\geq\frac{1}{4}|x-y|}\frac{C_1H(|u-z|^{-2})}{|u-z|^{d}}\leq \frac{C_14^dH(16|x-y|^{-2})}{|x-y|^{d}}
\leq \frac{C_14^{d+4}H(|x-y|^{-2})}{|x-y|^{d}}.\label{3.11}
\end{align}
By the same argument in (\ref {3.5}), we can apply Theorem \ref{mbhp} to get
\begin{equation}\label{3.15}
\P_x(X_{\tau_{U_1}} \in U_2) \leq c_{12}\P_{x_1}(X_{\tau_{U_1}} \in U_2)\frac{\delta_D(x)}{\delta_D(x_1)}\leq c_{13}\frac{\delta_D(x)}{\sqrt{t}}.
\end{equation}
 Applying (\ref{3.6}), (\ref{3.13}), (\ref{3.11}) and (\ref{3.15}) in the inequality (\ref{3.2}), we obtain
\begin{align*}
&p_D(t,x,y)\\
&\leq c_{14}\frac{\delta_D(x)}{\sqrt{t}}\left(t^{-d/2}\wedge\Big(t^{-d/2}e^{-\frac{|x-y|^2}{4C_2t}} + t\frac{H(|x-y|^{-2})}{|x-y|^{d}}+\phi^{-1}(t^{-1})^{d/2}e^{-\frac{a_U}{4}|x-y|^2\phi^{-1}(t^{-1})}\Big)\right)\\
&\quad +c_{15}\sqrt{t}\delta_D(x)\frac{H(|x-y|^{-2})}{|x-y|^{d}}\\
&\leq c_{16}\frac{\delta_D(x)}{\sqrt{t}} \left(t^{-d/2}\wedge\Big(t^{-d/2}e^{-\frac{|x-y|^2}{4C_2t}} + t\frac{H(|x-y|^{-2})}{|x-y|^{d}}+\phi^{-1}(t^{-1})^{d/2}e^{-\frac{a_U}{4}|x-y|^2\phi^{-1}(t^{-1})}\Big)\right).
\end{align*}
In the last line, (\ref{3.14}) is used. Combining   above two cases,  we have completed the proof of the proposition. 
\qed\\

{\bf Proof of Theorem \ref{T:3} (1) and Theorem \ref{T:4} (1).} We only prove Theorem \ref{T:3} (1), because
both proofs are almost identical.  Using (\ref{3.1}) and (\ref{e:upro}) instead of (\ref{3.2}) and (\ref{uu}) respectively, 
 we will follow  the proof of Proposition \ref{P:3.2}.
 
 Fix $T>0$. Let $t \in (0,T]$ and $x,y \in D$. By Proposition \ref{P:3.2}, Theorem \ref{T:2} and symmetry, we only need to prove Theorem \ref{T:3} (1) when $\delta_D(x) \vee \delta_D(y) < r_0\sqrt{t}/(16\sqrt{T}) \leq r_0/16$. Thus we assume that  $\delta_D(x) \vee \delta_D(y) < r_0\sqrt{t}/(16\sqrt{T}) \leq r_0/16$.
 Define $x_0, x_1$ and $U_1$ in the same way as in the proof of Proposition \ref{P:3.2} and let $c_1:=((d+1)/2) \vee ((dC_L^{-1/\gamma}T^{1/\gamma-1}\phi^{-1}(1)^{-1})/(C_2a_U))\vee (r_0^2/(16C_2T))$ where $a_U$ and $C_2$ are constants in \eqref{e:upro}. Now we estimate $p_D(t,x,y)$ by considering the following two cases separately.\\
{\it Case (1):} $|x-y| \leq 4(C_2c_1)^{1/2}\sqrt{t}$. By the semigroup property and symmetry, 
$$
p_D(t,x,y)=\int_Dp_D(t/2,x,z)p_D(t/2,z,y)dz\leq \Big(\sup_{z\in D}p_D(t/2,y,z)\Big)\int_Dp_D(t/2,x,z)dz.
$$
Using Proposition \ref{P:3.2} and \eqref{3.7} in the above inequality , we have
\begin{align}\label{de}
p_D(t,x,y) &\leq c_2t^{-d/2}\bigg(1\wedge \frac{\delta_D(y)}{\sqrt{t}}\bigg)\P_x(\tau_D>t/2)\nonumber\\
&\leq c_2t^{-d/2}\bigg(1\wedge \frac{\delta_D(y)}{\sqrt{t}}\bigg)\bigg(1\wedge \frac{\delta_D(x)}{\sqrt{t}}\bigg)\nonumber\\
&\leq c_2e^{c_1}t^{-d/2}e^{-|x-y|^2/(16C_2t)}\bigg(1\wedge \frac{\delta_D(x)}{\sqrt{t}}\bigg)\bigg(1\wedge \frac{\delta_D(y)}{\sqrt{t}}\bigg).
\end{align}
{\it Case (2):} $|x-y| \geq 4(C_2c_1)^{1/2}\sqrt{t}$. Define $U_2$ and $U_3$ in the same way as in (\ref{3.9}). Then by the same way  \eqref{3.12} and \eqref{3.11} hold. 

By our choice of $c_1$, we have $t\leq |x-y|^2/(8(d+1)C_2)$. Using this and the fact that $s \to s^{-(d+1)/2}e^{-\beta/s}$ is  increasing on the interval $(0,2\beta/(d+1)]$, we get for $s\leq t$,
$$
s^{-(d+1)/2}e^{-\frac{|x-y|^2}{16C_2s}}\leq t^{-(d+1)/2}e^{-\frac{|x-y|^2}{16C_2t}}.
$$ 
Thus by \eqref{e:upro} and \eqref{3.12},
\begin{align}\label{ue1}
&\sup_{s\leq t, z\in U_2}p_D(s,z,y)\nn\\
&\leq c_3\sup_{s\leq t, z\in U_2} 
 \bigg(1 \wedge \frac{\delta_D(y)}{\sqrt{s}}\bigg)\nn \\&\quad\times \left(s^{-d/2}\wedge \Big(s^{-d/2}e^{-\frac{|z-y|^2}{4C_2s}} + \frac{sH(|z-y|^{-2})}{|z-y|^{d}}+\phi^{-1}(s^{-1})^{d/2}e^{-\frac{a_U}{4}|z-y|^2\phi^{-1}(s^{-1})}\Big)\right)\nn\\
&\leq c_3\delta_D(y)\sup_{\substack{s\leq t\\ |z-y|\geq \frac{|x-y|}{2}}}\Big(s^{-(d+1)/2}e^{-\frac{|z-y|^2}{4C_2s}}+\frac{\sqrt{s}H(|z-y|^{-2})}{|z-y|^{d}}+\frac{\phi^{-1}(s^{-1})^{d/2}}{\sqrt{s}}e^{-\frac{a_U}{4}|z-y|^2\phi^{-1}(s^{-1})} \Big)\nn\\
&\leq c_4\delta_D(y)\sup_{s\leq t}\Big(s^{-(d+1)/2}e^{-\frac{|x-y|^2}{16C_2s}}+\frac{2^d\sqrt{t}H(4|x-y|^{-2})}{|x-y|^{d}}+s^{-1/2}\phi^{-1}(s^{-1})^{d/2}e^{-\frac{a_U}{16}|x-y|^2\phi^{-1}(s^{-1})} \Big)\nn\\
&\leq c_5\frac{\delta_D(y)}{\sqrt{t}}\Big(t^{-d/2}e^{-\frac{|x-y|^2}{16C_2t}}+\frac{tH(|x-y|^{-2})}{|x-y|^{d}}\Big) \nn\\
&\quad+c_4\delta_D(y)\Big(\sup_{s\leq t}s^{-1/2}e^{-\frac{a_U}{32}|x-y|^2\phi^{-1}(s^{-1})}\Big)
\Big(\sup_{s\leq t}
\phi^{-1}(s^{-1})^{d/2}e^{-\frac{a_U}{32}|x-y|^2\phi^{-1}(s^{-1})}\Big).
\end{align}
Now we find upper bounds of  $s^{-1/2}e^{-\frac{a_U}{32}|x-y|^2\phi^{-1}(s^{-1})}$ and
$
\phi^{-1}(s^{-1})^{d/2}e^{-\frac{a_U}{32}|x-y|^2\phi^{-1}(s^{-1})}$ for $s \le t$. By our choice of $c_1$ and \eqref{phisc}, we have 
$$t\leq \frac{a_U}{16d}C_L^{1/\gamma}T^{1-1/\gamma}\phi^{-1}(1)|x-y|^2\quad \text{and} \quad\phi^{-1}(t^{-1})^{-1}\leq C_L^{-1/\gamma}T^{1/\gamma-1}t\phi^{-1}(1)^{-1}\leq \frac{a_U}{16d}|x-y|^2.$$
 Using this and the fact that $s \to s^{-d/2}e^{-\beta/s}$ is increasing on the interval $(0,2\beta/d]$,  we get
\begin{align}\label{ue2}
&\delta_D(y)\big(\sup_{s\leq t}s^{-1/2}e^{-\frac{a_U}{32}|x-y|^2\phi^{-1}(s^{-1})}\big)
\big(\sup_{s\leq t}\phi^{-1}(s^{-1})^{d/2}e^{-\frac{a_U}{32}|x-y|^2\phi^{-1}(s^{-1})}\big)\nn\\
&\leq \delta_D(y)\big(\sup_{s\leq t}s^{-1/2}e^{-\frac{a_U}{32}|x-y|^2C_L^{1/\gamma}T^{1-1/\gamma}\phi^{-1}(1)s^{-1}}\big)
\big(\sup_{s\leq t}\phi^{-1}(s^{-1})^{d/2}e^{-\frac{a_U}{32}|x-y|^2\phi^{-1}(s^{-1})}\big)\nn\\
&\leq \delta_D(y)t^{-1/2}e^{-\frac{a_U}{32}C_L^{1/\gamma}T^{1-1/\gamma}\phi^{-1}(1)t^{-1}}\phi^{-1}(t^{-1})^{d/2}e^{-\frac{a_U}{32}|x-y|^2\phi^{-1}(t^{-1})}\nn\\
&= \frac{\delta_D(y)}{\sqrt{t}}e^{-\frac{a_U}{32}C_L^{1/\gamma}T^{-1/\gamma}\phi^{-1}(1)}\phi^{-1}(t^{-1})^{d/2}e^{-\frac{a_U}{32}|x-y|^2\phi^{-1}(t^{-1})}.
\end{align}
Therefore combine \eqref{ue1} with \eqref{ue2} we get
\begin{align}\label{3.18}
&\sup_{s\leq t, z\in U_2}p_D(s,z,y)\nn\\
&\leq c_6\frac{\delta_D(y)}{\sqrt{t}}\Big(t^{-d/2}e^{-\frac{|x-y|^2}{16C_2t}}+\frac{tH(|x-y|^{-2})}{|x-y|^{d}}+\phi^{-1}(t^{-1})^{d/2}e^{-\frac{a_U}{32}|x-y|^2\phi^{-1}(t^{-1})}\Big)\nn\\
&\leq c_7\frac{\delta_D(y)}{\sqrt{t}}\left(t^{-d/2}\wedge \Big(t^{-d/2}e^{-\frac{|x-y|^2}{16C_2t}}+\frac{tH(|x-y|^{-2})}{|x-y|^{d}}+\phi^{-1}(t^{-1})^{d/2}e^{-\frac{a_U}{32}|x-y|^2\phi^{-1}(t^{-1})}\Big)\right).
\end{align}
In the last inequality we have used similar argument as the one leading  (\ref{3.14}). On the other hand by (\ref{3.7}) we have
\begin{align}
&\int_0^t\P_x(\tau_{U_1}>s)\P_y(\tau_D>t-s)ds\leq \int_0^t\P_x(\tau_D>s)\P_y(\tau_D>t-s)ds\nonumber\\
&\leq c_8\int_0^t\frac{\delta_D(x)}{\sqrt{s}}\frac{\delta_D(y)}{\sqrt{t-s}}ds=c_8\delta_D(x)\delta_D(y)\int_0^1\frac{1}{\sqrt{r(1-r)}}dr=c_9\delta_D(x)\delta_D(y). \label{3.19}
\end{align}
Using (\ref{3.11}), (\ref{3.15}), (\ref{3.18}) and (\ref{3.19}) in the inequality (\ref{3.1}), we conclude that 
\begin{align*}
&p_D(t,x,y)\\
&\leq c_{10}\frac{\delta_D(x)\delta_D(y)}{t}\left(t^{-d/2}\wedge \Big(t^{-d/2}e^{-\frac{|x-y|^2}{16C_2t}}+\frac{tH(|x-y|^{-2})}{|x-y|^{d}}+\phi^{-1}(t^{-1})^{d/2}e^{-\frac{a_U}{32}|x-y|^2\phi^{-1}(t^{-1})}\Big)\right) \\
&\quad+c_{11}\frac{\delta_D(x)\delta_D(y)}{t}\times\frac{tH(|x-y|^{-2})}{|x-y|^{d}}\\
&\leq c_{12}\frac{\delta_D(x)\delta_D(y)}{t}\left(t^{-d/2}\wedge \Big(t^{-d/2}e^{-\frac{|x-y|^2}{16C_2t}}+\frac{tH(|x-y|^{-2})}{|x-y|^{d}}+\phi^{-1}(t^{-1})^{d/2}e^{-\frac{a_U}{32}|x-y|^2\phi^{-1}(t^{-1})}\Big)\right)\\
&=c_{12}\bigg(1\wedge \frac{\delta_D(x)}{\sqrt{t}}\bigg)\bigg(1\wedge \frac{\delta_D(y)}{\sqrt{t}}\bigg)\\
&\quad\times\left(t^{-d/2}\wedge \Big(t^{-d/2}e^{-\frac{|x-y|^2}{16C_2t}}+\frac{tH(|x-y|^{-2})}{|x-y|^{d}}+\phi^{-1}(t^{-1})^{d/2}e^{-\frac{a_U}{32}|x-y|^2\phi^{-1}(t^{-1})}\Big)\right).
\end{align*}
In the second inequality above, we also used similar argument as the one leading  (\ref{3.14}). This combined with (\ref{de}) completes the proof.
\qed\\
{\bf Proof of Theorem \ref{T:3} (3).}
The proof is same as \cite[Theorem 1.1(iii), (iv)] {CKS7}. We should consider Theorem \ref{T:3}(1) and Theorem \ref{T:3}(2) instead of \cite[Theorem 1.1(ii)]{CKS7} and \cite[Theorem 2.6]{CKS7}, respectively. We omit the proof.
\qed
\section{Green function estimates}
In this section we give the proof of Corollaries \ref{ge} and \ref{dge}. \\
{\bf Proof of Corollary \ref{ge}.}
Note that by Lemma \ref{lem:bf}, if $|x|\leq 1$, $\phi(1)\phi(|x|^{-2})^{-1} \geq |x|^2$, and if $|x|>1$, $\phi(1)\phi(|x|^{-2})^{-1}<|x|^2$. We split the integral
$$
G(x)=\int_{0}^{\phi(1)\phi(|x|^{-2})^{-1}\wedge |x|^2}p(t,x)dt +\int_{\phi(1)\phi(|x|^{-2})^{-1}\wedge |x|^2}^{\infty}p(t,x)dt.
$$
By Remark \ref{R:L2.3} and \eqref{e:weakupper},
\begin{align*}
&\int_{0}^{\phi(1)\phi(|x|^{-2})^{-1}\wedge |x|^2}p(t,x)dt \\
&\leq c_1\int_{0}^{\phi(1)\phi(|x|^{-2})^{-1}\wedge |x|^2} \Big(t^{-d/2}e^{-c_2|x|^2/t} +t|x|^{-d}H(|x|^{-2})+\phi^{-1}(t^{-1})^{d/2}e^{-c_3|x|^2\phi^{-1}(t^{-1})} \Big) dt \\
&\leq c_4\Bigg(\int_{0}^{\phi(1)\phi(|x|^{-2})^{-1}\wedge |x|^2} t^{-d/2}e^{-c_2|x|^2/t} dt +\int_{0}^{\phi(1)\phi(|x|^{-2})^{-1}\wedge |x|^2} t|x|^{-d}\phi(|x|^{-2}) dt \Bigg)\\
&\leq c_5\Bigg(\int_{0}^{\phi(1)\phi(|x|^{-2})^{-1}\wedge |x|^2} t^{-d/2}(|x|^2/t)^{-d/2} dt +|x|^{-d}\phi(|x|^{-2}) \int_{0}^{\phi(1)\phi(|x|^{-2})^{-1}\wedge |x|^2} tdt  \Bigg)\\
&\leq c_6\Big(|x|^{-d}\big(\phi(|x|^{-2})^{-1}\wedge |x|^2\big) + |x|^{-d}\big(\phi(|x|^{-2})^{-1}\wedge |x|^{4}\phi(|x|^{-2})\big)\Big)\\
&\leq c_7\Big(|x|^{-d+2}\wedge |x|^{-d}\phi(|x|^{-2})^{-1}\Big).
\end{align*}
For $|x|\leq 1$, using Lemma \ref{L:2.4},
\begin{align*}
\int_{\phi(1)\phi(|x|^{-2})^{-1}\wedge |x|^2}^{\infty}p(t,x)dt=\int_{ |x|^2}^{\infty}p(t,x)dt \leq c_8\int_{ |x|^2}^{\infty} t^{-d/2}dt=\frac{2c_8}{d-2}|x|^{-d+2}.
\end{align*}
For $|x|>1$, using Lemma \ref{L:2.3} and change of variables, 	
\begin{align*}
\int_{\phi(1)\phi(|x|^{-2})^{-1}\wedge |x|^2}^{\infty}p(t,x)dt &=\int_{\phi(1)\phi(|x|^{-2})^{-1}}^{\infty}p(t,x)dt\leq c_9\int_{\phi(1)\phi(|x|^{-2})^{-1}}^{\infty}\phi^{-1}(t^{-1})^{d/2}dt\\
&=c_9\int_0^{c_{10}|x|^{-2}}s^{d/2}\bigg(-\frac{1}{\phi(s)}\bigg)'ds=c_9\int_0^{c_{10}|x|^{-2}}s^{d/2-1}\frac{y\phi'(s)}{\phi(s)^2}ds\\
&\leq c_9\int_0^{c_{10}|x|^{-2}}\frac{s^{d/2-1}}{\phi(s)}ds,
\end{align*}
in the last inequality we use $\lambda\phi'(\lambda) \leq \phi(\lambda)$ because $\phi$ is represented by \eqref{lapex}.
Since $d\geq 3$, we have $\frac{d}{2}-2>-1$. Hence
$$
\int_0^{c_{10}|x|^{-2}}\frac{s^{d/2-1}}{\phi(s)}ds\leq \frac{c_{11}}{\phi(|x|^{-2})}\int_0^{c_{10}|x|^{-2}} s^{d/2-1}\bigg(\frac{|x|^{-2}}{s}\bigg)ds=c_{12}|x|^{-d}\phi(|x|^{-2})^{-1}.
$$
On the other hand, for $|x|>1$, by Lemma \ref{L:2.11} and the condition of $L_0(\gamma, C_L)$ on $\phi$, we have
\begin{align*}
G(x) &\geq \int_{\phi(|x|^{-2})^{-1}}^{\infty}p(t,x)dt \\
&\geq c_{13}\int_{\phi(|x|^{-2})^{-1}}^{2\phi(|x|^{-2})^{-1}}\phi^{-1}(t^{-1})^{d/2}dt \geq c_{13}\phi^{-1}\Big(\frac{\phi(|x|^{-2})}{2}\Big)^{d/2}\frac{1}{\phi(|x|^{-2})} \geq \frac{c_{13}(C_L/2)^{d/(2\gamma)}}{|x|^d\phi(|x|^{-2})}.
\end{align*}
When $|x|\leq 1$, by Lemma \ref{L:2.11} and Lemma \ref{L:2.12}, we have
\begin{align*}
G(x)&\geq \int_{|x|^2}^{2|x|^2}p(t,x)dt \geq c_{14}\int_{|x|^2}^{2|x|^2} t^{-d/2}\wedge \phi^{-1}(t^{-1})^{d/2}dt\geq c_{15}\int_{|x|^2}^{2|x|^2}t^{-d/2}dt \geq c_{15} \frac{|x|^{-d}}{2^{d/2}}|x|^2.
\end{align*}
Third inequality holds because for $t\leq 2$, $t^{-d/2} \leq c\phi^{-1}(t^{-1})^{d/2}$  for some $c_{16}>0$. Hence we conclude that
$$
G(x) \asymp |x|^{-d+2}\wedge |x|^{-d}\phi(|x|^{-2})^{-1}.
$$
\qed\\
{\bf Proof of Corollary \ref{dge}.} 
Recall that $g_D(x,y)$ is defined in Corollary \ref{dge}. Let $T:=$diam$(D)^2$. Then we have (see the proof of \cite[Corollary 1.3]{CKS7})
$$
\int_0^T \bigg(1 \wedge \frac{\delta_D(x)}{\sqrt{t}}\bigg)\bigg(1 \wedge \frac{\delta_D(y)}{\sqrt{t}}\bigg) p^{(2)}(t,c(x-y))dt+\int_T^{\infty} e^{-\lambda_1t}\delta_D(x)\delta_D(y)dt \asymp g_D(x,y).
$$
By Theorem \ref{T:3}(2) and (3),  we can easily obtain $G_D(x,y)\geq c_1g_D(x,y)$.

Next we consider the upper bound for $G_D(x,y)$.  By Theorem \ref{T:3} (1), for the bounded $C^{1,1}$-open set $D$,
\begin{align*}
&G_D(x,y)=\int_0^{\infty}p^D(t,x,y)dt\\
&\quad\leq c_1\int_0^T \bigg(1 \wedge \frac{\delta_D(x)}{\sqrt{t}}\bigg)\bigg(1 \wedge \frac{\delta_D(y)}{\sqrt{t}}\bigg)  \\
&\quad\quad  \times \left(t^{-d/2}\wedge \Big(p^{(2)}(t,c_2(x-y)) + \frac{tH(|x-y|^{-2})}{|x-y|^{d}}+\phi^{-1}(t^{-1})^{d/2}e^{-a_U|x-y|^2\phi^{-1}(t^{-1})}\Big)\right)\\
&\quad\quad +c_1\int_T^{\infty} e^{-\lambda_1t}\delta_D(x)\delta_D(y)dt\\
&\quad\leq  c_1\int_0^T \bigg(1 \wedge \frac{\delta_D(x)}{\sqrt{t}}\bigg)\bigg(1 \wedge \frac{\delta_D(y)}{\sqrt{t}}\bigg)  \\
&\quad \quad  \times \left(p^{(2)}(t,c_2(x-y))+\Big(t^{-d/2}\wedge \big(\frac{tH(|x-y|^{-2})}{|x-y|^{d}}+\phi^{-1}(t^{-1})^{d/2}e^{-a_U|x-y|^2\phi^{-1}(t^{-1})}\big)\Big)\right)\\
&\quad\quad +c_1\int_T^{\infty} e^{-\lambda_1t}\delta_D(x)\delta_D(y)dt\\
&\quad\leq c_2\left(g_D(x,y)+\int_0^T \bigg(1 \wedge \frac{\delta_D(x)}{\sqrt{t}}\bigg)\bigg(1 \wedge \frac{\delta_D(y)}{\sqrt{t}}\bigg) \big(t^{-d/2} \wedge \frac{t}{|x-y|^{d+2}} \big)dt\right).
\end{align*}
For the last inequality, we use \eqref{e:weakupper}, boundedness of $D$ and Lemma \ref{lem:bf}:$$
\frac{tH(|x-y|^{-2})}{|x-y|^{d}}+\phi^{-1}(t^{-1})^{d/2}e^{-a_U|x-y|^2\phi^{-1}(t^{-1})}\leq \frac{c_3t\phi(|x-y|^{-2})}{|x-y|^{d}} \leq \frac{c_4t}{|x-y|^{d+2}}.
$$
To complete the proof, it suffices to show that
$$
\int_0^T \bigg(1 \wedge \frac{\delta_D(x)}{\sqrt{t}}\bigg)\bigg(1 \wedge \frac{\delta_D(y)}{\sqrt{t}}\bigg) \big(t^{-d/2} \wedge \frac{t}{|x-y|^{d+2}} \le c_5 g_D(x,y),
$$
which is \cite[(4.1)]{CKS7}. Thus the remaining proof is same as the part of proof starting on the page 135 in \cite[Corollary 1.3]{CKS7}. So we omit it.
\qed

\bigskip
\noindent

\vspace{.1in}



\small
\begin{thebibliography}{99}

\bibitem{B} J.~Bertoin, L\'evy Processes, Cambridge University Press, Cambridge (1996).

\bibitem{BGR} K.~Bogdan, T.~Grzywny, M.~Ryznar, Barriers, exit time and survival probability for unimodal L\'evy processes. Probab. Theory Relat. Fields {\bf 162} (2015) 155-198.
\bibitem{CKK} Z.-Q.~Chen, P.~Kim, T.~Kumagai, Weighted Poincar\'e Inequality and Heat Kernel Estimates for Finite Range Jump Processes. Math. Ann. {\bf 342(4)} (2008), 833-883.


\bibitem{CKK2} Z.-Q.~Chen, P.~Kim, T.~Kumagai, Global heat kernel estimates for symmetric jump processes. Trans. Amer. Math. Soc. {\bf 363} (2011), 5021-5055.



\bibitem{CKS2} Z.-Q.~Chen, P.~Kim, R.~Song, Heat kernel estimates for Dirichlet fractional Laplacian. J. European Math. Soc., {\bf 12} (2010), 1307-1329.

\bibitem{CKS3} Z.-Q.~Chen, P.~Kim, R.~Song, Two-sided heat kernel estimates for censored stable-like processes. Probab. Theory Relat. Fields, {\bf 146} (2010), 361-399.

\bibitem{CKS4} Z.-Q.~Chen, P.~Kim, R.~Song, Dirichlet heat kernel estimates for $\Delta^{\alpha/2}+\Delta^{\beta/2}$. Ill. J. Math. {\bf 54} (2010), 1357-1392.

\bibitem{CK2} Z.-Q.~Chen, T.~Kumagai, Heat kernel estimates for stable-like processes on d-sets. Stochastic Process. Appl. {\bf 108} (2003), 27-62.

\bibitem{CK3} Z.-Q.~Chen, T.~Kumagai, Heat kernel estimates for jump processes of mixed types on metric measure spaces. Probab. Theory Relat. Fields, {\bf 140} (2008), 277-317.



\bibitem{CK} Z.-Q.~Chen, T.~Kumagai, A priori H\"older estimate, parabolic Harnack principle and heat kernel estimates for diffusions with jumps. Rev. Mat. Iberoamericana {\bf 26} (2010), 551-589

\bibitem{CKS5} Z.-Q.~Chen, P.~Kim, R.~Song, Heat kernel estimate for $\Delta + \Delta^{\alpha/2}$ in $C^{1,1}$ open sets. J. London Math. Soc. {\bf 84(1)}  (2011), 58-80.

\bibitem{CKS6} Z.-Q.~Chen, P.~Kim, R.~Song, Sharp heat kernel estimates for relativistic stable processes in open sets. Ann. Probab. {\bf 40} (2012), 213-244.

\bibitem{CKS7} Z.-Q.~Chen, P.~Kim, R.~Song, Dirichlet Heat Kernel Estimates for Subordinate Brownian Motion with Gaussian Components. {J. Reine Angew. Math. } {\bf 711} (2016) 111--138.

\bibitem{Ch} S.~Cho, Two-sided global estimates of the Green's function of parabolic equations. {Potential Analysis} {\bf 25(4)} (2006) 387--398.

\bibitem{CZ} K.L.~Chung, Z.~Zhao, From Brownian Motion to Schr\"odinger's Equation, Springer, Berlin, 1995.
\bibitem{KM} P.~Kim, A.~Mimica, Estimates of Dirichlet Heat Kernel for Subordinate Brownian Motions.  	arXiv:1708.08606 [math.PR]

\bibitem{KSV2} P.~Kim, R.~Song, Z.~Vondra\v{c}ek, On the potential theory of one-dimensional subordinate Brownian motions with continuous components. {Potential Anal.} {\bf 33} (2010), 153-173.


\bibitem{KSV1} P.~Kim, R.~Song, Z.~Vondra\v{c}ek, Potential theory of subordinate Brownian motions with Gaussian components. {Stoch.~Processes Appl.} {\bf 123(3)} (2013) 764--795.



\bibitem{M}
A.~Mimica, Heat kernel estimates for subordinate brownian motions. Proc. London Math. Soc., {\bf 113(5)} (2016) 627-648.




\bibitem{Sat} K.-i.~Sato, L\'evy processes and infinitely divisible distributions, Cambridge University Press, 1999.




\bibitem{SV}
R.~Song, Z.~Vondra\v{c}ek, Parabolic Harnack inequality for the mixture of Brownian motion and stable process. {Tohoku Math. J. (2)} {\bf 59(1)} (2007), 1--19.

\bibitem{Zh} Q.~S.~Zhang, The boundary behavior of heat kernels of Dirichlet Laplacians. {J. Differential Equations} {\bf 182} (2002) 416--430.

\end{thebibliography}
\end{document}